\documentclass[]{interact}

\usepackage{epstopdf}
\usepackage{subfigure}
\usepackage{graphicx}
\usepackage{mathtools}
\usepackage{amssymb}
\usepackage{amsthm}
\usepackage{thmtools}
\usepackage{xcolor}
\usepackage{nameref}
\usepackage{bm}
\usepackage{hyperref}
\usepackage{url}
\usepackage{subcaption}

\usepackage{natbib}
\bibpunct[, ]{(}{)}{;}{a}{}{,}
\renewcommand\bibfont{\fontsize{10}{12}\selectfont}

\theoremstyle{plain}

\theoremstyle{definition}

\theoremstyle{remark}
\newtheorem{remark}{Remark}

\usepackage{subcaption}
\renewcommand{\r}{{\bm r}}
\newcommand{\Om}{{\bm \Omega}}
\newcommand{\R}{\mathbb{R}}
\renewcommand{\S}{\mathbb{S}}
\newcommand{\del}{\partial}
\newcommand{\n}{{\bm n}}
\renewcommand{\t}{{\bm t}}
\renewcommand{\b}{{\bm b}}
\renewcommand{\v}{{\bm v}}
\newcommand{\h}{{\bm h}}
\renewcommand{\u}{{\bm u}}
\newcommand{\V}{{\bm V}}

\newcommand{\kappan}{{\bm{\kappa}^{\n}}}
\newcommand{\kappat}{{\bm{\kappa}^{\t}}}
\newcommand{\kappab}{{\bm{\kappa}^{\b}}}
\newcommand{\kappabt}{{\bm{\kappa}^{\tilde{\b}}}}
\newcommand{\curl}{\text{rot}}
\newcommand{\X}{{\bm X}}
\newcommand{\N}{{\bm N}}
\newcommand{\T}{{\bm T}}
\let\oldkappa\kappa
\renewcommand{\kappa}{\bm{\oldkappa}}

\let\oldgamma\gamma
\renewcommand{\gamma}{\bm{\oldgamma}}
\newcommand{\sinphi}{\sin(\phi)}
\newcommand{\sinomega}{\sin(\omega)}
\newcommand{\sintheta}{\sin(\theta)}
\newcommand{\cosphi}{\cos(\phi)}
\newcommand{\cosomega}{\cos(\omega)}
\newcommand{\costheta}{\cos(\theta)}
\newcommand{\quaddd}{\quad\quad\quad}

\begin{document}

\articletype{Preprint}

\title{Curvilinear coordinates and curvature in radiative transport \\ \small (Submitted for review to Journal of computational and theoretical transport on 8/28/25)}

\author{
\name{Johannes Krotz\textsuperscript{a}\thanks{Johannes Krotz, Email: jkrotz@nd.edu} and Ryan G. McClarren\textsuperscript{a}}
\affil{\textsuperscript{a} Department of Aerospace and Mechanical Engineering, University of Notre Dame, Notre Dame, 46545, IN, USA  }
}

\maketitle

\begin{abstract}
We derive a general expression for the streaming term in radiative transport equations and other transport problems when formulated in curvilinear coordinates, emphasizing coordinate systems adapted to the geometry of the domain and the directional dependence of particle transport. By parametrizing the angular variable using a local orthonormal frame, we express directional derivatives in terms of curvature-related quantities that reflect the geometry of underlying spatial manifolds. Our formulation highlights how the interaction between coordinate choices and curvature influences the streaming operator, offering geometric interpretations of its components. The resulting framework offers intuitive insight into when and how angular dependence can be simplified and may guide the selection of coordinate systems that balance analytical tractability and computational efficiency.
\end{abstract}


\section{Introduction} \label{sec:Introduction}

	Equations of the form  
\begin{align}
	\frac{1}{c}\frac{\partial \Psi}{\partial t} + \boldsymbol{\Omega} \cdot \nabla \Psi + \mathcal{G}[\Psi] = Q, \label{eq:kinetic_general}
\end{align}
arise in a wide variety of physical and engineering contexts involving the directional transport of particles, energy or other agents. Here, $\Psi = \Psi(\mathbf{r}, \boldsymbol{\Omega}, t)$ denotes a density or distribution function depending on spatial position $\mathbf{r} \in \mathbb{R}^3$, direction of movement or momentum $\boldsymbol{\Omega} \in \mathbb{S}^2$, and time $t$. The streaming term $\boldsymbol{\Omega} \cdot \nabla \Psi$ captures advective transport along characteristics determined by $\boldsymbol{\Omega}$, while the operator $\mathcal{G}[\Psi]$ accounts for local or nonlocal interactions such as absorption, scattering, collisions, or relaxation effects. The source term $Q$ models external or internal production mechanisms.

Equation~\eqref{eq:kinetic_general} and closely related formulations find application in a wide range of domains. Examples range from  rarefied gas dynamics where it appears as the Boltzmann equation describing non-equilibrium molecular flows~\cite{cercignani1988}, to plasma physics, where the Vlasov equation models the evolution of charged particle distributions under electromagnetic fields~\cite{nicholson1983}. The phonon Boltzmann equation describes microscale heat transport in solids~\cite{majumdar1993}, and similar forms govern the propagation of cosmic rays in astrophysics~\cite{schlickeiser2002}. Radiative transfer models photon transport in stellar and atmospheric environments~\cite{chandrasekhar1960}, as well as in optical tomography~\cite{arridge1999}. While in the context of biological and social systems, kinetic equations with directional transport can be used to describe traffic flow, swarming, and crowd dynamics~\cite{bellomo2013}.

Despite their disparate origins, these systems all share a common mathematical core structured around the advective transport term $\boldsymbol{\Omega} \cdot \nabla \Psi$. This term fundamentally couples direction and space, giving rise to anisotropy, long-range effects, and complex interactions that distinguish kinetic equations from classical PDEs in purely spatial domains.

	As with many problems in physics and engineering, it is common practice to address this problem using coordinates tailored to the geometry of the domain. For instance, spherical coordinates for scenarios where most parameters exhibit rotational symmetry. Notably, equation \eqref{eq:kinetic_general} offers additional degrees of freedom in choosing coordinates for $\Om$, beyond those for $\r$. However, this introduces added complexity: the optimal representation of $\r$ and $\Om$ may depend on each other.
	
\begin{figure}[h]
		\centering
		\includegraphics[width=0.5\linewidth,trim={8cm 10cm 5cm 4.5cm},clip]{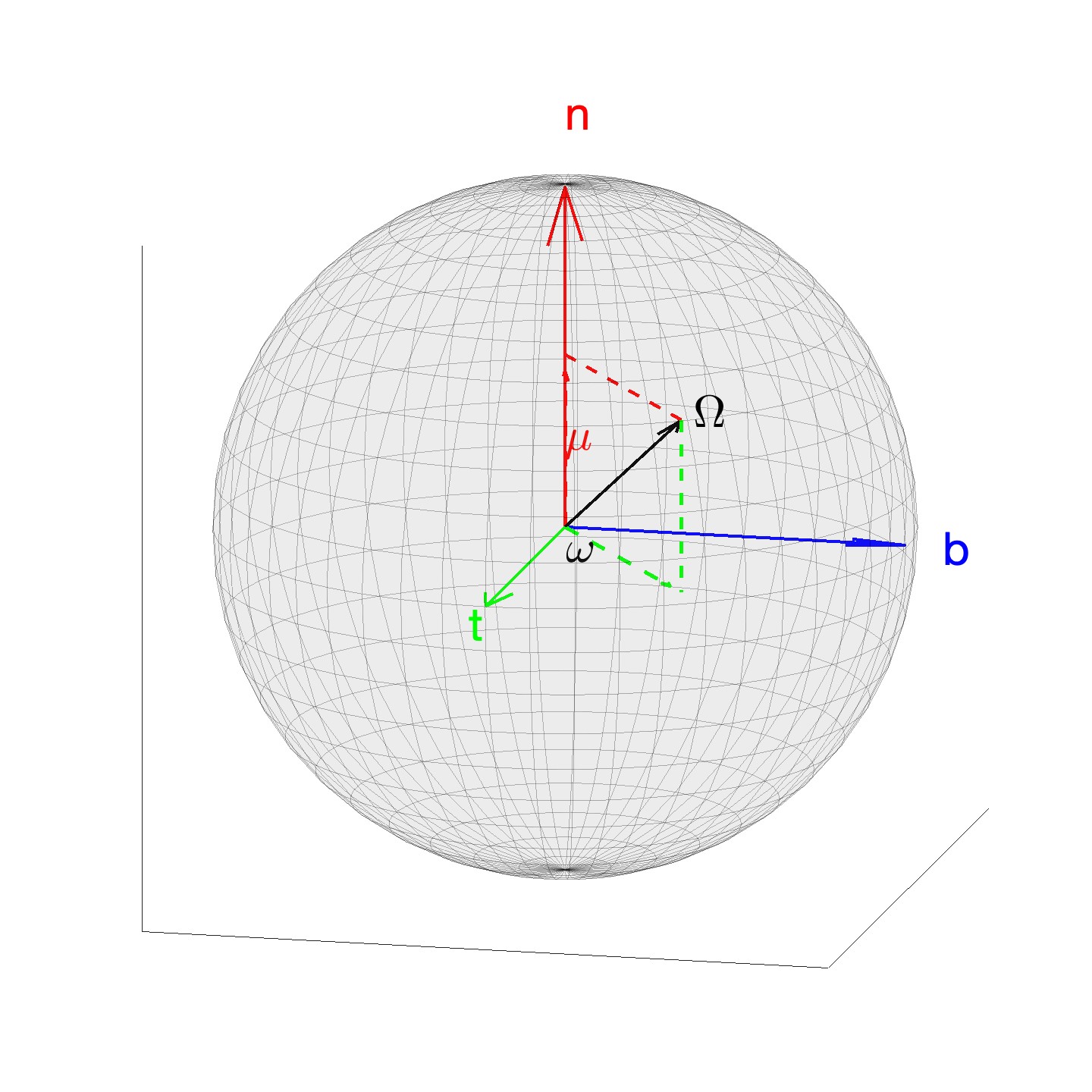}
		\caption{Visualization of the parametrization of $\Om(\mu,\omega)\in\S^2$: $\mu=\Om\cdot \n$ is the projection of $\Om$ onto $\n$, $\omega$ measures angle from $\t$ in $\t-\b$-plane. $\n,\t,\b$ might depend on $\r$ making $\mu$ amd $\omega$ functions of $\r$ even for fixed $\Om$  }
		\label{fig:1}
	\end{figure}

	For $\Om$, a vector on the unit sphere, a common parametrization is the standard spherical coordinates (polar and azimuthal angles) or equivalently:
	\begin{align}
		\Om(\omega,\mu) = \mu  \n+  \sqrt{1-\mu^2}\cos(\omega) \t \nonumber \\\quaddd+ \sqrt{1-\mu^2}\sin(\omega)\b , \label{eq:Om}
	\end{align}
	where $\omega \in [0,2\pi]$, $\mu \in [-1,1]$, and $\n, \t, \b \in \S^2$ form an orthogonal basis. While this formulation may appear convoluted, it simply means that $\Om$ is parametrized by its projection $\mu$ onto a designated vector $\n$ and a polar angle $\omega$ measured in a plane orthogonal to $\n$ from a reference vector $\t$ in that plane. A visualization is provided in Figure~\ref{fig:1}. 
	
	Given a parametrization, i.e., a set of coordinates for $\r \in \R^3$, say $\r = \r(u^1,u^2,u^3)$, a natural choice for $\n, \t$, and $\b$ would be the canonical tangent vectors $\n=c_1\frac{\del \r}{\del{u^1}}$, $\t = c_2\frac{\del \r}{\del u^2}$, and $\b =\n\times \t$, where the coefficients $c_{1},c_2$ ensure normalization. Note that this canonical choice is not necessarily orthogonal; if it is not, additional modifications are required. Orthogonalization via Gram-Schmidt is one possible approach. Other choices in terms of contra and covariant coordinates of the canonical tangent vectors, guaranteeing orthogonality, are described in \cite{FREIMANIS}. 
	Moving forward, we will not rely on any specific choice of $\n, \t$, and $\b$, but will assume they form an integrable orthonormal vector field.
	
	The parametrization of $\Om$ presented in \eqref{eq:Om} simplifies integrals over $\Om$, which frequently appear in scattering operators $\mathcal{S}$, as:
	\begin{align}
		\int_{\S^2}(\cdot) \, d\Om = \int_{-1}^1 \int_0^{2\pi}(\cdot) \, d\omega \, d\mu.
	\end{align}
	A more significant advantage, however, is that $\n$ can often be chosen as the normal vector to level sets of $\Psi$, leading to cases where $\Psi$ is effectively independent of $\omega$, i.e.,
	\begin{align}
		\Psi(\r,\Om,E,t) = \Psi(\r,(\omega,\mu),E,t) = \Psi(\r,\mu,E,t).
	\end{align}
	This reduction in the dimensionality of the parameter space can greatly simplify solving equations \eqref{eq:kinetic_general}. However, achieving this simplification with constant $\n$ is generally only possible in trivial cases. In most scenarios, $\n, \t$, and $\b$ must depend on $\r$, making the computation of the term $\Om \cdot \nabla \Psi$ particularly cumbersome.
	
	In the literature, the exact form of this term for a given choice of coordinates for $\r$ and $\Om$ is often derived through tedious geometric arguments. These derivations can be challenging to follow and even harder to adapt to slightly different geometries. Examples of such derivations can be seen in \cite{Bell_Glasstone,helical}
	
	In the next section, we derive general expressions for $\Om \cdot \nabla \Psi$ for arbitrary choices of $\r$ and $\n, \t, \b$ in a purely analytical manner. 
	We denote by $\nabla^{\r}$ the partial derivatives in space and identify 
	\begin{align}
		\Om \cdot \nabla \Psi = \Om\cdot \nabla^{\r} \Psi + \frac{\del \Psi}{\del \mu}\Om\cdot \nabla \mu  + \frac{\del \Psi}{\del \omega}\Om\cdot \nabla \omega, \label{eq: OmPsi}
	\end{align}
	presenting alternative formulations for the terms $\Om \cdot \nabla \mu$ and $\Om \cdot \nabla \omega$ in section \ref{sec:Derivation_main}. 
    We relate these expressions to the geometric properties of curves and surfaces embedded in $\R^3$, which we believe will make the process of deriving $\Om \cdot \nabla \Psi$ more straightforward, in some cases, and provides additional intuition and guidance in general.
    Namely we will rewrite the terms $\Om\cdot \nabla \mu$ and $\Om\cdot \nabla \omega$ as terms looking as follows:\begin{align}
    \Om\cdot \nabla \mu = (1-\mu^2) \fbox{
  \begin{tabular}{c}
    \text{The curvature of} \\
      \text{some manifold} $X$
  \end{tabular}} + \mu\sqrt{1-\mu^2} \fbox{  \begin{tabular}{c}
  \text{The curvature of}\\ \text{a curve orthogonal to $X$}
  \end{tabular}}
    \end{align}
    and 
\begin{align}
\Om \cdot \nabla \omega = \sqrt{1-\mu^2} \fbox{
  \begin{tabular}{c}
    \text{The curvature of} \\
      \text{some curves in} $X$
  \end{tabular}}  + \mu \fbox{
  \begin{tabular}{c}
    \text{Non-curvature term} \\
      \text{quantifying spinning of $\t$ around $\n$} 
  \end{tabular}}.
\end{align}
    
     In \cite{FREIMANIS} equation \eqref{eq: OmPsi} is calculated for arbitrary coordinates $\r$, but restricted to two choices of $\n,\t,\b$ relative to the chosen coordinates, and with focus on the coordinate representation alone. We pose severely fewer restrictions on $\n,\t,\b$ and focus on interpretations regarding the embeddings into $\R^3$ rather than purely on the coordinate level. 
	In section \ref{sec:examples}, we provide examples illustrating this approach in some classical settings and some less common scenarios.

	\section{Derivation of ${\bf \Omega}\cdot \nabla \Psi$}\label{sec:Derivation_main}
	
	In this section, we set out to calculate an expression $\Om\cdot\nabla \Psi$ when $\Om$ is parametrized as 
	\begin{align}
		\Om(\r) =\mu \n+ \sqrt{1-\mu^2}\cos(\omega) \t + \sqrt{1-\mu^2}\sin(\omega)\b,
	\end{align}
	where $\n,\t\in C^1(\R^3,\R^3)$, and $\b(\r):=\n\times \t$ are orthonormal, differentiable vector fields and thus $\Psi(\r,\Om,E,t) = \Psi(\r,(\omega(\r),\mu(\r)),E,t)$.
	Given coordinates for $\r$ it is not always the simplest choice to guarantee all these conditions. Seemingly obvious choices might have neither normality nor orthogonality. We will stick to these assumptions regardless and discuss possible fixes, if they do not hold in some of the examples in section \ref{sec:examples}.
	
	We will begin by establishing some notation and simple properties of these vector fields, as well as $\omega$ and $\mu$, that will be helpful in later calculations. We will then proceed to calculate $\Om\cdot\nabla \Psi$, providing several alternative forms for most terms involved, allowing practitioners to choose the one best suited to their needs.
	Throughout the bold letters will be used for vectors in $\R^3$. 	Further, given a vector $\h_1,\h_2\in \R^3$, we denote by $\h_1\cdot \h_2$ their standard inner product, by $||\h_1|| = \sqrt{\h_1\cdot \h_1}$ the euclidean norm and  their  normalization as $\hat{\h_1}:= \frac{\h_1}{||\h_1||}$.We will denote the directional derivative in direction $\h_1$ as $\nabla_{\h_1}:= \h_1\cdot \nabla$. 
	Lastly, let us define the projection of $\Om(\r)$ onto the orthogonal complement of $\n$, i.e., onto the plane spanned by $\t$ and $\b$, as $\Om_{\parallel}(\r)$. Specifically, 
	\begin{align}
		\Om_{\parallel}(\r):= & \sqrt{1-\mu^2}\cos(\omega) \t + \sqrt{1-\mu^2}\sin(\omega)\b 
		 =  \Om(\r)-\mu\n.
	\end{align}

		We will analyze equation \eqref{eq: OmPsi} piece by piece and start with the term $\Om\cdot\nabla\mu = \nabla_{\Om}\mu$ in section \ref{sec:nabla_mu}.  
        
        Between $\nabla_{\Om}\mu$ and $\nabla_{\Om}\omega$ the prior is arguably more important, since terms stemming from the second one can often be ignored if $\n$ is chosen orthogonal to level sets of $\Psi$, which under mild conditions leads to $\del\Psi/\del \omega=0.$ (Compare \ref{app:1})
	For cases where  $\del\Psi/\del \omega\neq 0$ we will derive an expression for $\nabla_{\Om}\omega$ in section \ref{sec:nabla_omega}. Most terms in the streaming 	will be related to curvatures of either curves or surfaces, or both, defined by the vectors $\n,\t$ and $\b$, which lends itself to a more intuitive approach to choosing coordinates.

	\subsection{Derivation of $ {\bf\Omega}\cdot\nabla\mu$}\label{sec:nabla_mu}
In this section we set out to further evaluate and interpret the term $\nabla_{\Om}\mu$ in equation \eqref{eq: OmPsi}.
	
	Since $\mu = \Om\cdot \n$ we first find \begin{align}
	&	\nabla_{\Om} \mu  = \Om\cdot \nabla_{\Om} \n  = \Om_\parallel \cdot (\nabla_{\Om_{\parallel}}\n  + \mu\nabla_{\n}\n) + \mu \n\cdot \nabla_{\Om}\n\label{eq:musplit} 
	\end{align}
	using the linearity of the derivative and the inner product.  Note that the last term $\mu \n\cdot \nabla_{\Om}\n =0,$ because $\n$ is assumed to be normalized. 
	It also is a consequence of a more general property of orthonormal frames detailed in the following remark.

	\begin{remark}{{Properties of orthonormal frames}}\label{rem:orth_frames}\\
		Let $\h \in \S^2$ and let $\u:\R^3\to\S^2$ and $\v:\R^3\to\S^2$ be distinct, differentiable orthonormal vector fields, i.e.
		$\u\cdot \u=\v\cdot \v=1$ and $\u\cdot \v = 0$ everywhere. Since these inner products are constant in space their derivatives disappear and therefore
		\begin{align}
			\u \cdot \nabla_{\h}\u =0
		\end{align}
		and
		\begin{align}
			\u\cdot\nabla_{\h} \v  = - \v\cdot\nabla_{\h} \u.
		\end{align}  
		Note that this holds for $\u,\v\in\{\n,\t,\b\}.$
	\end{remark}

	Using Remark \ref{rem:orth_frames} we are left to simplify
	\begin{align}
		\nabla_{\Om} \mu  = \Om_\parallel \cdot \nabla_{\Om_{\parallel}}\n  + \mu(\Om_{\parallel} \cdot \nabla_{\n}\n)  \label{eq:12}
	\end{align}
	The second term can be expressed in terms of the normal curvature of an integral curve.
	\begin{remark}{Integral Curves and Normal Curvature}\\\label{rem:curv_curvature}
		Given a differentiable vector field $\u: \R^3\to\S^2$ let $\gamma^{\u}(\r,\tau)$ be the integral curve of $\u$ passing through $\r$,  i.e. let it solve the following differential equation:
		\begin{align}
			\frac{d \gamma^{\u}(\tau,\r)}{d\tau} &= \u(\gamma^{\u}(\tau,\r))  & \quad \tau\in (-\varepsilon,\varepsilon) \nonumber\\
			\gamma^{\u}(0,\r) & = \r   &
		\end{align}
		Existence of $\gamma^{\u}$ is guaranteed by Picard-Lindel\"off, but since we will only use this curve locally we simply note that we can solve this differential equation in a neighborhood of $\gamma(0,\r)$. We shall further, sloppily, refer to $\gamma^{\u}(0,\r)$ as $\gamma^{\u}(\r)$.  Since $\u$ is a unit vector everywhere $\gamma^{\u}$ is parametrized by arc length and has a normal curvature vector 
		\begin{align}
			\kappa^{\u}(\r)&:= -\frac{d^2 \gamma^{\u}(\tau,\r)}{d \tau^2}\bigg|_{\tau=0}  =- \frac{		d \u(\gamma^{\u}(\tau,\r))}{d\tau}\bigg|_{\tau=0} =-(\u(\r)\cdot \nabla) \u(\r) = - \nabla_{\u}\u(\r).
		\end{align}
		(The sign convention here is chosen such that were $\gamma^{\u}$ a curve on a sphere with outward normal $\bm N$, the value of $\kappa^{\u}\cdot \bm N$ would be positive.)
	\end{remark}
	Setting $\u = \n$ in the previous remark, we find that 	the second term in equation \eqref{eq:12} is
	
	\begin{align}
		\mu	(\Om_\parallel \cdot \nabla_{\n} \n ) &=  \mu\sqrt{1-\mu^2}(\hat{\Om}_{\parallel}\cdot \nabla_{\n}\n) = \mu\sqrt{1-\mu^2}\left(\hat{\Om}_{\parallel} \cdot \bm{\kappa}^{\n}\right)
	\end{align}
	expressing the term entirely in terms of $\mu$ and the curvature of the integral curve of $\n$ in the direction of $\Om_{\parallel}$. 
	This term quantifies how much $\gamma^\n$ curves in direction $\Om_{\parallel}$. It disappears, if $\gamma^\n$ is a straight line. If we imagine a set of surfaces orthogonal to $\gamma^\n$ along the curve, this is simultaneously a measure for how parallel such surfaces are along the curve.
	
\begin{figure}\centering
    \subfigure[Spherical Coordinates]{
    \resizebox*{.45\textwidth}{!}{   
			\includegraphics[width=.45	\textwidth,trim={3cm 3cm 0cm 3cm},clip]{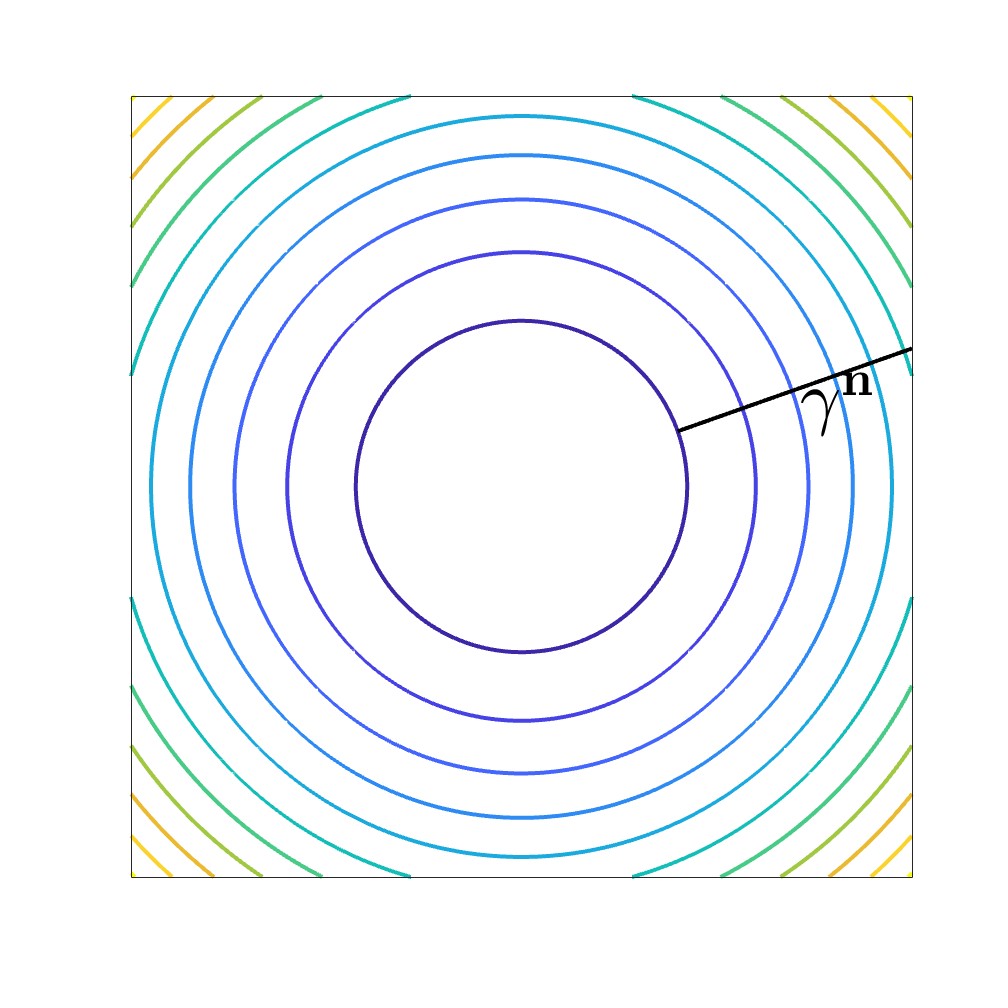}}}
            \subfigure[Elliptical Coordinates]{
    \resizebox*{.45\textwidth}{!}{
    
			\includegraphics[width=.45	\textwidth,trim={3cm 3cm 0cm 3cm},clip]{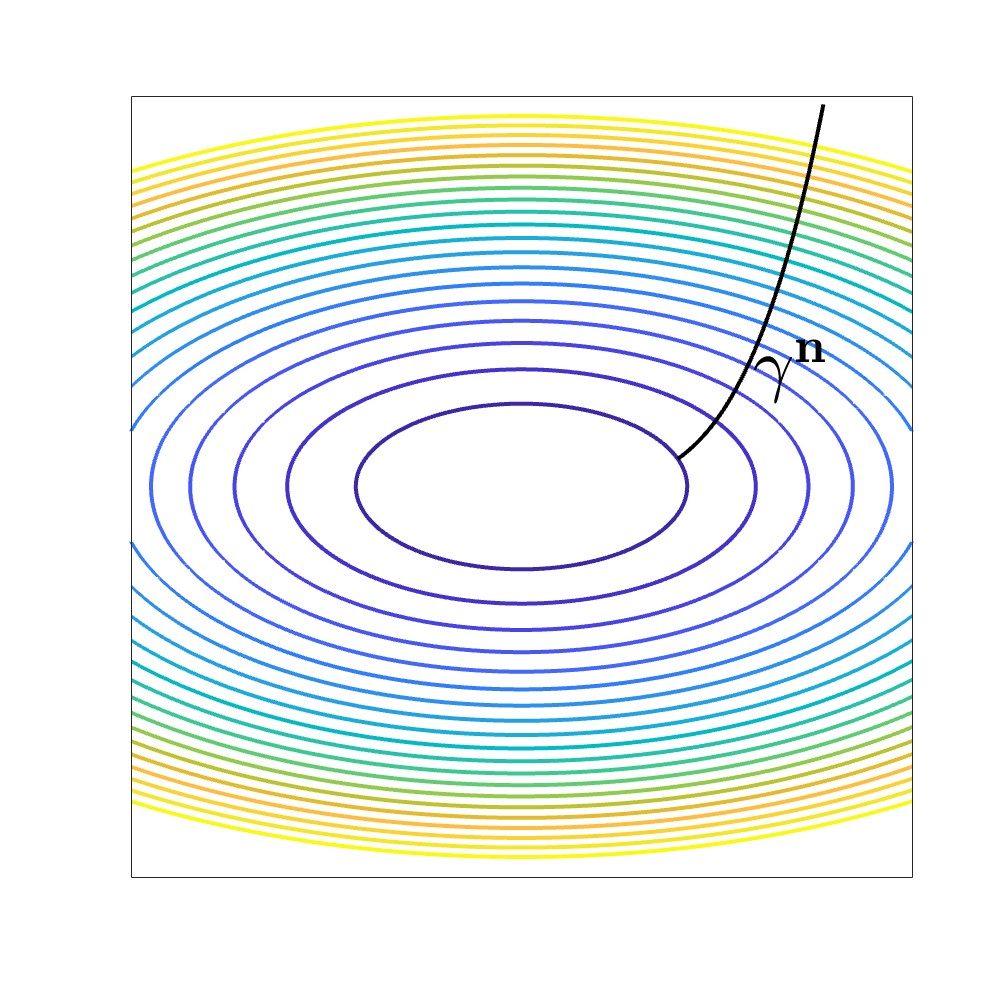}}}
		\caption{2D visualization of $\gamma^\n$ in spherical coordinates and elliptical coordinates. The integral curves $\gamma^\n$ following the normal vectors of circles are straight lines, while those of ellipses are curved.}	\label{fig:gamma_n}
\end{figure}

	As an example, if $\n(\r)$ is chosen to be the normalized gradient of $f(\r)=||\r||^2$, the term $\Om_{\parallel}\cdot \kappan = 0$ everywhere, since the level sets of $f$ are spheres centered at the origin and all these spheres are, locally parallel along $\gamma^{\n}$. On the other hand, if $\n$ is chosen to be the normalized gradient of $f(\r)= \r^\bot {\bf A}\r$, where ${\bf A}$ is some positive definite matrix with at least two distinct eigenvalues, then $\Om_{\parallel}\cdot\kappan$ will not vanish. In the latter case, the level sets of $f$ are ellipses, excluding spheres, which are not parallel along $\gamma^{\n}$. This is illustrated in figure \ref{fig:gamma_n}. We will see these two specific examples in section \ref{sec:examples}.
	
	Having established the notion of curvature of integral curves we shall first express the remaining term in equation \eqref{eq: OmPsi} in terms of integral curves as well. 	Recall the definition of $\Om_{\parallel}$ and therefor note that its normalization is given by $\hat{\Om}_{\parallel}=\cos(\omega)\t+\sin(\omega)\b$. 
	Thus using the linearity of the derivative and the inner product again we can now also tackle the first term in equation \eqref{eq:12} and see that
	
	\begin{align}
	&	\Om_\parallel \cdot \nabla_{\Om_{\parallel}} \n = (1-\mu^2)\left(\hat{\Om}_\parallel\cdot \nabla_{\hat{\Om}_{\parallel}} \n \right)  \nonumber\\
	&	=(1-\mu^2)\left[ \cos(\omega)^2(\t\cdot \nabla_{\t}\n)+ \sin(\omega)^2(\b\cdot\nabla_{\b}\n)  +\sin(\omega)\cos(\omega)\left( \t \cdot \nabla_{\b} \n + \b\cdot \nabla_{\t}\n \right) \right]  \end{align}
	Using Remark \ref{rem:orth_frames} the last set of derivatives can be rewritten as $\t \cdot \nabla_{\b} \n + \b\cdot \nabla_{\t}\n =  \t \cdot \nabla_{\b} \n - \t \cdot \nabla_{\b} \n = 0,$ causing the last summand to disappear. Further using Remark \ref{rem:orth_frames} again we can replace $\t\cdot \nabla_{\t}\n = -\n \cdot \nabla_{\t}\t$  and $\b\cdot \nabla_{\b}\n = - \n\cdot \nabla_{\b}\b$ and get
	
	\begin{align}
		&\Om_\parallel \cdot \nabla_{\Om_{\parallel}} \n \nonumber\\&= -(1-\mu^2)\left[ \cos(\omega)^2(\n\cdot \nabla_{\t}\t)+ \sin(\omega)^2(\n\cdot\nabla_{\b}\b)\right] \nonumber \\
		&=(1-\mu^2)\left[ \cos(\omega)^2 \n\cdot \kappat + \sin(\omega)^2 \n \cdot\kappab  \right].  
	\end{align}
	In the last step Remark \ref{rem:curv_curvature} was used to identify $\nabla_{\t}\t=\kappat$ and $\nabla_{\b}\b=\kappab$ as the curvatures of the corresponding integral curves.  
	
	While the interpretation in terms of the curvature of integral curves already offers some intuition, we would like to note that $\n\cdot \kappat$ and $\n\cdot \kappab$ can also be interpreted as the curvature of a manifold orthogonal to $\n$ in direction $\t$ and $\b$ respectively, if such a manifold exists. We will continue to look at 	$\Om_\parallel \cdot \nabla_{\Om_{\parallel}} \n$ a little longer to make its relation to the shape of such manifolds orthogonal to $\n$ clearer. To that end, we will establish under which conditions such orthogonal manifolds exist and  introduce the shape operator of a manifold.

\begin{figure}
		\centering
		\includegraphics[width=0.5\linewidth, trim={5cm 4cm 5cm 2.5cm},clip]{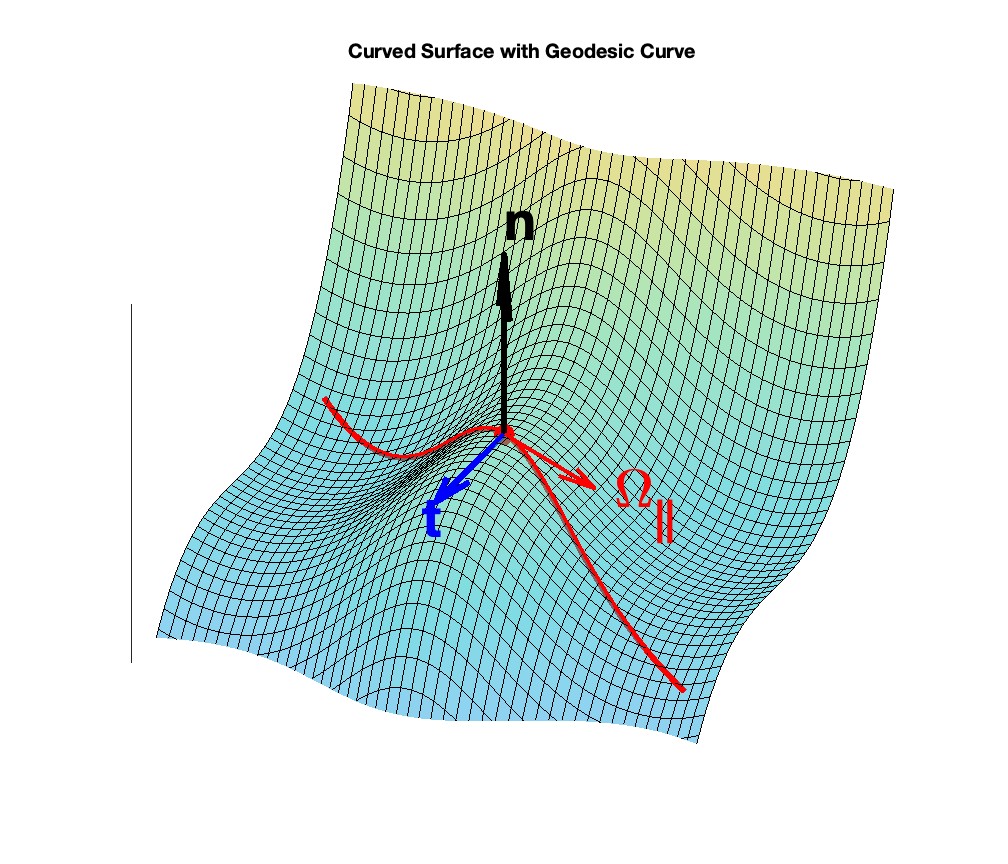}
		\caption{Example manifold $\X^\n$ with normal vector $\n$ and tangent vectors $\t$ and $\Om_{\parallel}$. The curvature at $\r$ in direction $\omega$ relative to $\t$, denoted by $C(\r,\omega)$ is the normal curvature of red curve starting at $\r$.}
		\label{fig:shape}
	\end{figure}

	\begin{remark}{Induced Foliation}\label{remark:Foliation}\\
		Let ${\V}:\R^3\to\S^2$ be a differentiable normalized vector field and let $U \subset \R^3$ be a connected open set. If ${\V}\cdot \curl ( \V)  = 0$ on $U$, then for each $\r\in U$ there is a manifold $X^\V_\r$ such that $\r \in X^\V_\r$ and so that $\V(\r')$ is the normal vector of $\X$ at all $\r' \in \X^\V_\r$. Since $\r' \in \X_\r^\V$ implies $\X_\r^\V= \X^\V_{\r'}$ this introduces an equivalence relation on $U$ called a foliation, with equivalence classes, called leaves, being surfaces normal to $\V$.  Due to this we will drop the $\r$ from $\X_\r^\V$ going forward and write $\X^\V$ to mean the leaf of the foliation that is orthogonal to $\V$ at $\r$. In $\R^3$ the condition $\V\cdot \curl(\V)=0$ is both necessary and sufficient for the existence of a foliation. This follows from Froebenius theorem \cite{Frobenius,Lee00} applied to vector fields in $\R^3$.

		A maybe more applicable, yet still sufficient (not necessary condition), for the existence of $X^\V$ can be expressed as follows: Let $\r(u^1,u^2,u^3)$ be a (local) parametrization of $\R^3$, i.e. let $(u^1,u^2,u^3)$ be the (local) coordinates. 
		If $0\neq\V\propto \frac{\del \r}{\del u^1}$ then  
		$\X^\V(u^2,u^3) = \r(\bar{u}^1,u^2,u^3)$ 
		is a parametrization of $\X^\V$ containing the point $\bar{\r}=\r(\bar{u}^1,\bar{u}^2,\bar{u}^3)$, i.e. fixing one of the coordinates of $\r$, but altering the others generates a surface containing $\r$. 
	\end{remark}
	\begin{remark}{Manifolds and their Shape Operator}\label{remark:3}\\
		Let $\X^\N:\R^2\supset U\rightarrow \R^3$ be the parametrization of an embedded surface $\X(u^1,u^2)$ with parameters $(u^1,u^2)\in U$ and tangent vectors $\T_{i} = \frac{\del \X}{\del u^i}$ ($i=1,2$) and normal vector $\N\in \R^3.$ For simplicity assume $\N$ is not just defined and differentiable on $\X$, but on an open neighborhood $U$ of of $\X\subset U \subset \R^3$.
		The metric tensor $g$ of $\X$ can be expressed in coordinates as the matrix $g_{\alpha,\beta} = \T_\alpha \cdot \T_\beta$ where $\alpha,\beta \in\{u^1,u^2\}$, while the second fundamental form $h$ of $\X$ is given by the matrix defined through $h_{\alpha,\beta} = -\frac{\del \T_\alpha}{\del \beta}\cdot \N = -\frac{\del^2 \X}{\del \alpha \del \beta}\cdot \N$.
		From $g$ and $h$ the shape operator, or Weingarten map $S_\X$, can be calculated as $S_\X = g^{-1}h$. The shape operator contains information about the curvature of the manifold $\X$. Given a unit tangent vector $\u\in T_\X$ the curvature in direction $\u$ can be calculated as $\u^\top S_\X \u$. The shape operator shares its Eigenvalues with $h$, meaning Eigenvalues of $\S_X$ are the principal curvatures of $\X$, its determinant is the Gauss curvature and its trace is the mean curvature of $\X$.\\
		By Weingarten's theorem, the coordinate representation of $S_X$ can also be calculated as \begin{align}
			(S_X)_\alpha^\beta = \T_\alpha \cdot \nabla_{\T_\beta}\N & & \alpha,\beta \in \{u^1,u^2\}. 
		\end{align}		\cite{Weingarten,kreyszig2013}
	\end{remark} 
	With notation from Remarks \ref{remark:Foliation} and \ref{remark:3} and $\N=\n$ and $\T_{1}=\t$ and $\T_{2}=\b$ we can interpret the derivative terms in equation \eqref{eq:musplit} as components of the shape operator $S_{\X^\n}$ of the manifold  $\X^\n$ and thus
	\begin{align}
		&\Om_\parallel \cdot \nabla_{\Om_{\parallel}} \n \nonumber \\& = (1-\mu^2)\left[(\cosomega, \sinomega)\cdot S_{\X^\n} \cdot \begin{pmatrix}
			\cosomega \\ \sinomega
		\end{pmatrix}\right] \nonumber\\&=: (1-\mu^2)C(\r,\omega)
	\end{align}
	Note that $(\cosomega, \sin(\omega))$ is merely $\hat{\Om}_{\parallel}$ restricted to the tangent space of $\X^\n$  and that therefore the term in brackets is nothing but the curvature of $\X^\n$ in direction $\hat{\Om}_{\parallel}$, which we denoted by $C(\r,\omega)$. A visualization can be seen in figure \ref{fig:shape}.
	
	Overall we found two curvature forms for $\nabla_{\Om}\mu$ in 
	\begin{align}
		\nonumber\nabla_{\Om} \mu &=   
		   (1-\mu^2)\left[\cos(\omega)^2\left(\n\cdot \kappat\right)+\sin(\omega)^2\left(\n\cdot \kappab\right)  \right] +\mu\sqrt{1-\mu^2}\left(\hat{\Om}_\parallel\cdot \kappan\right)  \\
		&= (1-\mu^2)C(\r,\omega)+\mu\sqrt{1-\mu^2}\left(\hat{\Om}_\parallel\cdot \kappan\right).   
	\end{align}
	We will refer to the first one of these as the curve-curvature form of $\nabla_{\Om}\mu$ and to the second one as the surface-curvature form of $\nabla_{\Om}\mu$ when using them in examples later on. 
	
	\subsection{Derivation of ${\bf\Omega}\cdot \nabla \omega$ }\label{sec:nabla_omega}
\begin{figure}
		\centering
		\includegraphics[width=.5\linewidth, trim={14cm 12cm 14cm 10cm},clip]{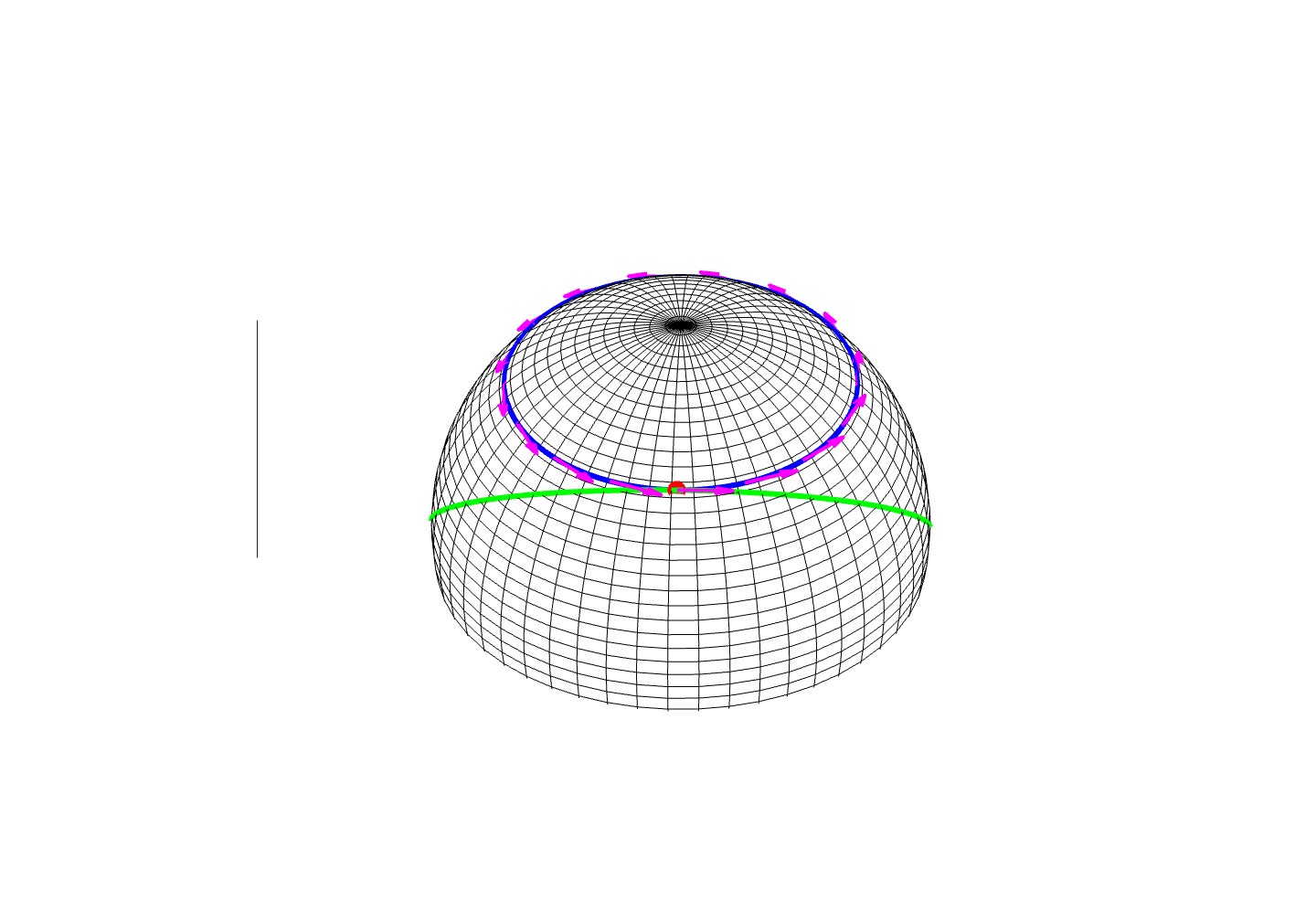}
		\caption{Example of non-geodesic integral curve $\gamma^\u$ (blue) along a tangent vector field $\u$ (magenta). For reference a geodesic starting at a common point (green). Since $\gamma^\u$ is non-geodesic the curvature $\kappa^u$ projected onto the tangent space is not zero.}
		\label{fig:sphere_geo}
	\end{figure}

	We will continue to calculate the remaining term in equation \eqref{eq: OmPsi}, namely $\nabla_{\Om}\omega$. To this end note that $\omega = \arccos\left(\frac{\Om_\parallel\cdot \t}{\sqrt{1-\mu^2}}\right) = \arccos\left(\hat{\Om}_{\parallel}\cdot \t\right)$ or alternatively $\omega = \arcsin\left(\frac{\Om_\parallel\cdot \b}{\sqrt{1-\mu^2}}\right) = \arcsin\left(\hat{\Om}_{\parallel}\cdot \b\right)$.  Between these two notions of $\omega$ $\omega$, the angle is uniquely determined up to multiples of $2\pi$. We will see that both lead to the same derivative.
	With $\arcsin(x)'=-\arccos(x)'=\frac{1}{\sqrt{1-x^2}}$ one sees either 
	\begin{align}
		\nabla_{\Om}\omega = \frac{1}{\cos(\omega)}\left(\hat{\Om}_\parallel\cdot \nabla_{\Om} \b \right) = \t\cdot \nabla_{\Om} \b \label{eq: omega1}
	\end{align}
	or \begin{align}
		\nabla_{\Om}\omega = \frac{1}{\sin(\omega)}\left(\hat{\Om}_\parallel\cdot \nabla_{\Om} \t \right) = -\b\cdot \nabla_{\Om}\t \label{eq:omega2}
	\end{align}
	In both of these equations, we used Remark \ref{rem:orth_frames}, specifically that $\b\nabla \b = \t\nabla \t = 0$, along with the definition of $\Om_{\parallel}$ to get to the right-hand side. Using Remark \ref{rem:orth_frames} again one can also conclude that equations \eqref{eq: omega1} and \eqref{eq:omega2} are indeed the same, where they are both defined.
To provide additional interpretation and intuition we will manipulite this further to also derive curvature forms for this equation.	
	First using the linearity of the derivative 
	\begin{align}
	&\nabla_{\Om}\omega =  \t\nabla_{\Om}\b = -\b\nabla_{\Om}\t \nonumber\\
	&	= \sqrt{1-\mu^2}\left[\cosomega (\t\cdot \nabla_{\t}\b) + \sinomega (\t\cdot \nabla_{\b}\b) \right] \nonumber + \mu (\t\cdot\nabla_{\n}\b) \label{eq:omega3}
	\end{align}
	is obtained. This is then rewritten, using Remark \ref{rem:orth_frames} and the definition of $\kappa^{\u}$ in Remark \ref{rem:curv_curvature}, as
	\begin{align}
		&\nabla_{\Om}\omega 
		= \sqrt{1-\mu^2}\left[-\cosomega (\b\cdot \nabla_{\t}\t) + \sinomega (\t\cdot \nabla_{\b}\b) \right] + \mu (\t\cdot\nabla_{\n}\b)\nonumber\\
		&=\sqrt{1-\mu^2}\left[\cosomega(\b\cdot \kappat)- \sinomega (\t\cdot \kappab)\right] \nonumber  +  \mu (\t\cdot\nabla_{\n}\b) \nonumber\\
		&= \frac{\del \Om_{\parallel}}{\del \omega}\cdot(\kappat+\kappab)+\mu(\t\cdot \nabla_\n \b)
	\end{align}
	The terms $\b\cdot \kappat$ and $\t\cdot\kappab$ are the curvatures of $\gamma^{\t}$ and $\gamma^{\b}$ within  $\X^\n$. They are zero if $\gamma^{\t}$ and $\gamma^{\b}$ are geodesics and otherwise measure how rapidly these curves bend away from geodesics. An example where $\gamma^\b$ is not a geodesic and such $\kappa^\b\cdot \t\neq 0$ is depicted in figure \ref{fig:sphere_geo}
	 
\begin{figure}
        \centering
        \includegraphics[width=.5\linewidth,trim={7cm 6cm 7cm 6cm},clip]{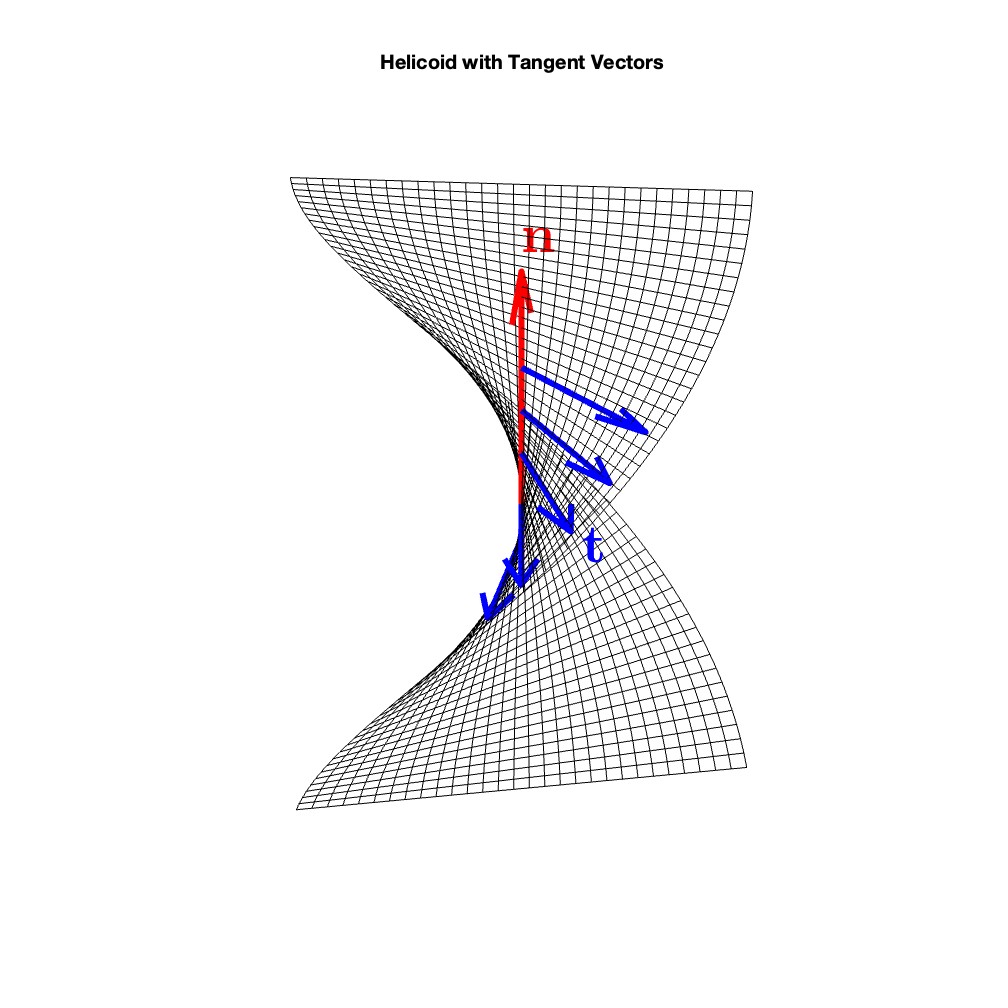}
        \caption{Example of a twisting, but flat manifold (helicoid). Diagonal entries of the shape operator disappear, despite off-diagonal terms being non-zero. The tangent vector $\t$ spins around $\n$ when moving along $\gamma^]n$.}
        \label{fig:helicoid}
    \end{figure}

	Notably the final term $\t\cdot\nabla_{\n}\b=- \b \cdot \nabla_{n}\t$ was not rewritten in terms of a curvature. It expresses how $\t$ and $\b$ spiral around the curve $\gamma^{\n}$. \\
	The previous section raises the question of whether a surface-curvature form of $\nabla_{\Om}$ exists as well.  Going back to equation \eqref{eq:omega3} we can indeed identify the terms $\t\cdot\nabla_\t \b$ and $\t\cdot \nabla_\n \b$ in terms of the shape operator $S_{\X^\b}$ of the surface $\X^\b$, if it exists.(Compare conditions of existence in Remark \ref{remark:Foliation}). The equation then becomes  
	\begin{align}
		\nonumber&\nabla_\Om \omega = (1,0) S_{\X^\b}\begin{pmatrix}
			\cosomega\sqrt{1-\mu^2}\\ \mu
		\end{pmatrix}  \sinomega\sqrt{1-\mu^2}(\t \cdot \nabla_\b\b)\\ & \nonumber=   (1,0) S_{\X^\b}\begin{pmatrix}
			\cosomega\sqrt{1-\mu^2}\\ \mu
		\end{pmatrix} -  \sinomega\sqrt{1-\mu^2}(\t\cdot \kappab) \label{eq:surface_b}
	\end{align}
	Alternatively rewriting $\t\cdot \nabla_\n\b = -\b\cdot \nabla_\n\t$ one can arive at a formulation in terms of $S_{\X^\t}$, again, if the surface $\X^\t$ exists.
	\begin{align}
		&\nabla_\Om \omega = - (1,0) S_{\X^\t}\begin{pmatrix}
			\sinomega\sqrt{1-\mu^2}\\ \mu
		\end{pmatrix}  - \cosomega\sqrt{1-\mu^2}(\t \cdot \nabla_\b\b)\\ & =   -(1,0) S_{\X^\t}\begin{pmatrix}
			\sinomega\sqrt{1-\mu^2}\\ \mu
		\end{pmatrix} +  \cosomega\sqrt{1-\mu^2}(\b\cdot \kappat)
	\end{align}
	Unfortunately, terms including the shape operator are not as symmetric as in the first section, so the interpretation becomes less straightforward. For example, while $(1,0)S_{\X^\b}(1,0)^\top$ still describes the curvature of $\X^\b$ in the first tangent direction, here in direction $\t$, the terms proportional to $(1,0)S_{\X^{\b}}(0,1)^\top$ correspond to a twisting of $\X^\b$ in direction of the second tangent vector, when moving along the first tangent vector. An example visualizing such twisting behavior is depicted in figure \ref{fig:helicoid}. Depicted is a scenario, where $\X^\b$ is a helicoid an intrinsically flat, but twisting surface, for which diagonal entries of $S_{\X^\b}$ would be zero, while the off-diagonal entries are not.

	\section{Examples}\label{sec:examples}
    In the previous section we found that \eqref{eq: OmPsi} can, in one of several ways, be recast as
    \begin{align}
		\Om \cdot \nabla \Psi &=\nonumber  
                \Om\cdot \nabla^{\r} \Psi \\&+ \left[(1-\mu^2)C(\r,\omega) + \mu\sqrt{1-\mu^2}\left( \hat{\Om}_\parallel\cdot\kappan\right) \right]\frac{\del \Psi}{\del \mu} \\&\nonumber+ \left[\frac{\del \Om_{\parallel}}{\del \omega}\cdot(\kappat+\kappab)+\mu(\t\cdot \nabla_\n \b)\right]\frac{\del \Psi}{\del \omega}. \label{eq:curv_form}
        \end{align}
        Before looking at specific examples, let us note that, for general curves and manifolds, curvatures can be rather complicated to calculate. However, the curvature approach lends itself nicely to a quick and intuitive assessment of which terms in \eqref{eq:curv_form} to expect in the first place. At first glance equation \eqref{eq:curv_form} should be considered as 
 \begin{align}
		\Om \cdot \nabla \Psi &=  
                \Om\cdot \nabla^{\r} \Psi \\&+ \left[(1-\mu^2) \fbox{
  \begin{tabular}{c}
    \text{Is} $X^\n$ \\
    \text{flat?}
  \end{tabular}} + \mu\sqrt{1-\mu^2}\fbox{
  \begin{tabular}{c}
    \text{Is} $\gamma^\n$ \text{curved} \\
      \text{towards} $X^\n$?
  \end{tabular}} \right]\frac{\del \Psi}{\del \mu} \\& + \left[
  \fbox{\begin{tabular}{c}
  \text{Are} $\gamma^\t$ \text{and} $\gamma^\b$ \\
  \text{curved in} $X^\n$? 
  \end{tabular}}
  +\mu   \fbox{\begin{tabular}{c}
  \text{Are} $\gamma^\t$ \text{and} $\gamma^\b$ \\
  \text{spinning around } $\gamma^\n$? 
  \end{tabular}}\right]\frac{\del \Psi}{\del \omega}. 
        \end{align}
If either of the boxed questions in the expression above is answered positively the corresponding term is known to be zero without any further calculation. Which question to ask for each term can be split further or altered based on whether a surface-curvature or a curve-curvature formulation is chosen in the previous section. Either way we believe this to be a quick and intuitive way to decide, if a term needs to be calculated explicitly or not.

    We will continue this section with a handful of examples, calculating all non-zero terms, but following this approach to skip the calculation of all terms that will end up being zero anyway.    
	Various choices for the parametrization of $\Om$ are explored to see how the introduced frameworks can make finding the coordinate form of $\Om\cdot \nabla \Psi$, easier and more intuitive. We will start with classical known examples and build up from there to increasingly more complicated scenarios. 
	\subsection{Cylindrical I}\label{sec:ex_cyl1}
\begin{figure}
        \centering
        \includegraphics[width=0.5\linewidth,trim={7cm 6cm 5cm 3cm},clip]{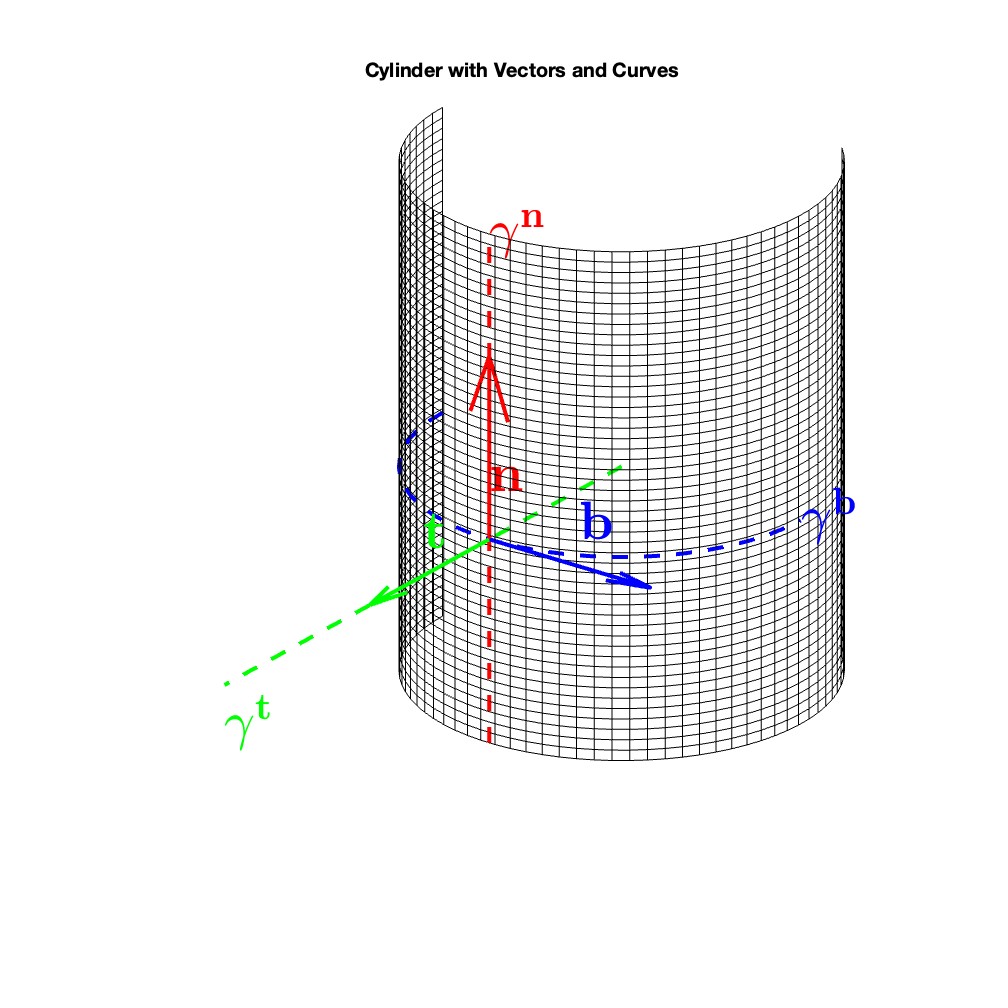}
        \caption{Visualization of $\Om$-basis choice $\n,\t,\b$ in cylindrical symmetry.}
        \label{fig:cylinder}
    \end{figure}

    Let us begin with a classical example useful in cases with cylindrical symmetry. Assume that $\r$ is expressed as $\r = (\rho \cos(\phi), \rho \sin(\phi),z)$ through the coordinates $\rho,\phi$ and $z$. Common choices in the literature are equivalent to picking $\n(\r) = (0,0,1)$ constant along the $z-$axis, $\t(\r) = (\cos(\phi), \sin(\phi),0)$ to be the outward normal vector of cylinders of radius $\rho$ and therefore $\b = (-\sin(\phi),\cos(\phi),0)$ the tangent vector of those cylinders in the $x-y-$plane. \cite{Bell_Glasstone} We show this choice in figure \ref{fig:cylinder}.
    Based on the set up we realize that a good choice to express equation \eqref{eq: OmPsi} here is as 
        \begin{align}
		\Om \cdot \nabla \Psi &=\nonumber  
                \Om\cdot \nabla^{\r} \Psi \\&+ \left[(1-\mu^2)C(\r,\omega) + \mu\sqrt{1-\mu^2}\left( \hat{\Om}_\parallel\cdot\kappan\right) \right]\frac{\del \Psi}{\del \mu} \\&\nonumber+ \left[   (1,0) S_{\X^\b}\begin{pmatrix}
			\cosomega\sqrt{1-\mu^2}\\ \mu
		\end{pmatrix} -  \sinomega\sqrt{1-\mu^2}(\t\cdot \kappab) \right]\frac{\del \Psi}{\del \omega},\label{eq:cyl_form}
        \end{align}
    which right away leads us towards the informal   

            \begin{align*}
		\Om \cdot \nabla \Psi &=\nonumber  
                \Om\cdot \nabla^{\r} \Psi \\&+ \left[(1-\mu^2)\fbox{
  \begin{tabular}{c}
    \text{Is} $X^\n$ \\
    \text{flat?}
  \end{tabular}} + \mu\sqrt{1-\mu^2}\fbox{
  \begin{tabular}{c}
    \text{Is} $\gamma^\n$ \text{curved} \\
      \text{towards} $X^\n$?
  \end{tabular}} \right]\frac{\del \Psi}{\del \mu} \\&\nonumber+ \left[   (1,0) \fbox{
  \begin{tabular}{c}
    \text{Is} $X^\b$ \\
    \text{flat?}
  \end{tabular}}\begin{pmatrix}
			\cosomega\sqrt{1-\mu^2}\\ \mu
		\end{pmatrix} -  \sinomega\sqrt{1-\mu^2}\fbox{
  \begin{tabular}{c}
    \text{Is} $\gamma^\b$ \text{curved} \\
    \text{towards $\t$?}
  \end{tabular}} \right]\frac{\del \Psi}{\del \omega}.
        \end{align*}
	The answer to the first three questions will be 'Yes', telling us that the corresponding terms will disappear without further calculation. Here is why. Since $\n$ is constant, it is straightforward to see, that all $\X^\n$ are planes and therefor flat.  Similarly, $\gamma^\n$ are straight lines in $z-$direction. Thus $\kappan=0$ as well, leading to $\nabla_\Om \mu =0$ every where.  \\
	 Since $\b$ is independent of both $\rho$ and $z$ the surface $\X^\b$ are just planes as well. The only term remaining is the fourth one, which in this case can also be determined without direct calculations.
     The curves $\gamma^\b$ generated by $\b$ are great circles around cylinders of radius $\rho$ with normal vector $\t$. Therefore $\t\cdot\kappab = \frac{1}{\rho}$ is given by the principal curvature of such cylinders. We therefore concluded $\nabla_{\Om}\omega = -\sinomega\sqrt{1-\mu^2}/\rho$ and hence
	
	\begin{align}
		\nabla_{\Om}\Psi = \nabla^\r_{\Om} \Psi+ 0 \frac{\del \Psi}{\del \mu} - \frac{\sqrt{1-\mu^2}\sin(\omega)}{\rho}\frac{\del \Psi}{\del \omega}. \end{align}
	Expressing $\nabla^\r_\Om \Psi$ in cylinder coordinates also makes sense in practice, but is not required. 
	\subsection{Cylindrical II}\label{sec:ex_cyl2}
	While the previous example is a common choice it is not the only way $\Om$ could be parametrized if the underlying system has cylindrical symmetry. Here is an example, where the roles of $\n,\t$ and $\b$ are switched compared to the previous example. Assume $\r$ again to be given as $\r = (\rho \cos(\phi), \rho \sin(\phi),z)$. Now let $\n$ be the radial vector $(\cosphi,\sinphi,0)$, i.e. the normal vector to cylinders of radius $\rho$. Let $\t = (-\cos(\phi),\sinphi,0)$ be the tangent vector to these cylinders orthogonal to the $z-$axis and let $\b=(0,0,1)$ be the parallel to the $z-$axis. Since $\b$ is constant, we can see immediately, that both $S_{\X^\b}$ and $\kappab$ are $0$. Therefore $\nabla_{\Om}\omega = 0$.\\
	The trade-off is that, since $\n$ is now the normal vector on a cylinder everywhere, i.e. $X^\n$ are cylinders, the corresponding shape operators do not disappear anymore. It is in fact given by 
	\begin{align}
		S_{X^\n} = \begin{pmatrix}
			1/\rho & 0 \\ 0 & 0
		\end{pmatrix}.
	\end{align}
	Since $\gamma^\n$ are still straight lines and thus $\kappan = 0$ we have $\nabla_{\Om}\mu = \sqrt{1-\mu^2}\cos (\omega)^2/\rho$ or \begin{align}
		\nabla_{\Om}\Psi = \nabla_{\Om}^\r\Psi + \frac{\sqrt{1-\mu^2}\cos(\omega)^2}{\rho} \frac{\del\Psi}{\del\mu} + 0\frac{\del \Psi}{\del \omega}.
	\end{align}
	
	\subsection{Sphere}\label{sec:ex_sphere}
\begin{figure}
        \centering
        \includegraphics[width=0.5\linewidth,trim={8cm 6cm 6cm 8cm},clip]{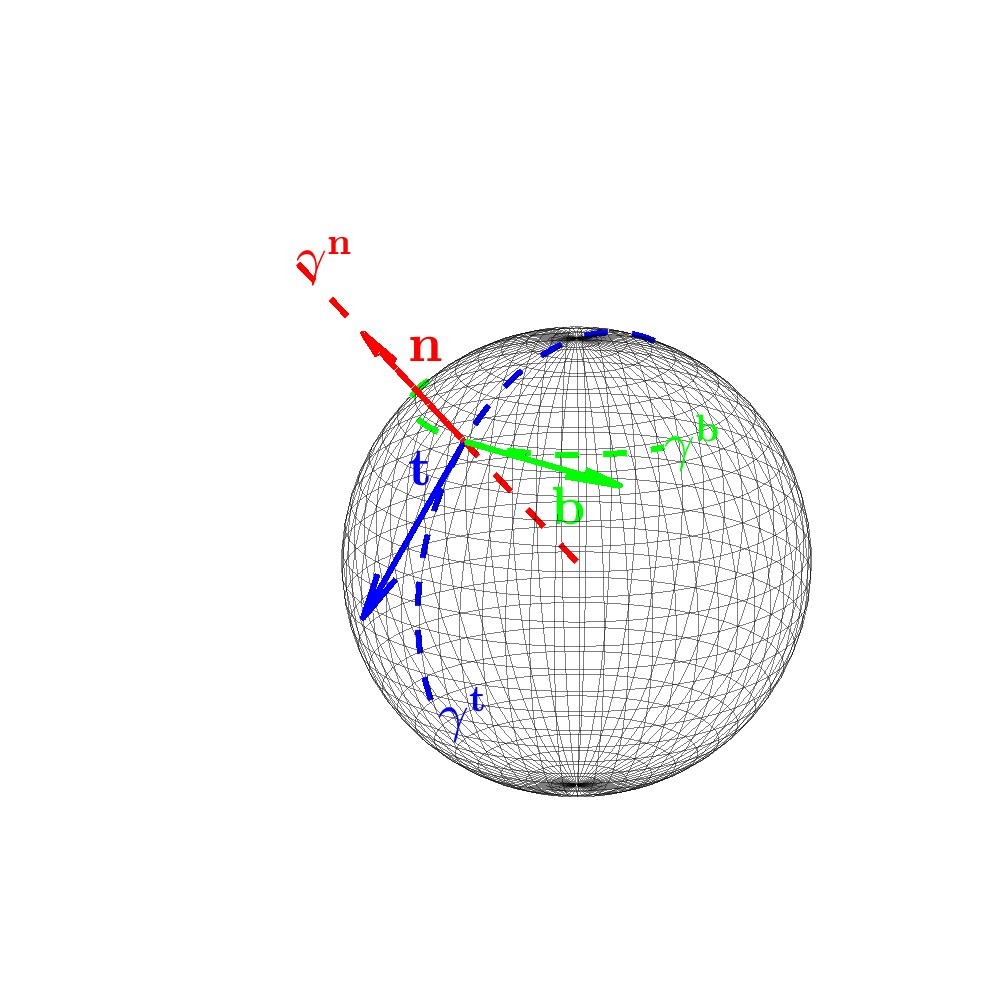}
        \caption{Visualization of $\Om$-basis choice $\n,\t,\b$ in spherical symmetry.}
        \label{fig:sphere_om}
    \end{figure}

	We will continue with another classical example: spherical coordinates. Assume 
	$\r = \rho (\cos(\phi)\sin(\theta), \sin(\phi)\sin(\theta),\cos(\theta))$. The vector $\n(\r)$ is commonly chosen to be the normal vector to the sphere of radius $\rho$ at $\r$, i.e.    
	$\n = (\cos(\phi)\sin(\theta), \sin(\phi)\sin(\theta),\cos(\theta)).$  The tangent vectors $\t$ and $\b$ are picked as $\t = (\cos(\phi)\cos(\theta),\sin(\phi)\cos(\theta),\sin(\theta))$, which is the tangent vector of the geodesic connecting north and south pole under the angle $\theta$ and $\b = (\cos(\phi),\sin(\phi),0)$, parallel to curves circling the poles parallel to the equator. The choice of $\n,\t$ and $\b$ is depicted in figure \ref{fig:sphere_om} for clarity.
    
    We express equation \eqref{eq: OmPsi} here is as 
        \begin{align}
		\Om \cdot \nabla \Psi &=\nonumber  
                \Om\cdot \nabla^{\r} \Psi \\&+ \left[(1-\mu^2)C(\r,\omega) + \mu\sqrt{1-\mu^2}\left( \hat{\Om}_\parallel\cdot\kappan\right) \right]\frac{\del \Psi}{\del \mu} \\&\nonumber+ \left[ \sqrt{1-\mu^2}\left[\cosomega(\b\cdot \kappat)- \sinomega (\t\cdot \kappab)\right] \nonumber  +  \mu (\t\cdot\nabla_{\n}\b) \right]\frac{\del \Psi}{\del \omega},\label{eq:cyl_form}
        \end{align}
which gives the informal version 
            \begin{align*}
		\Om \cdot \nabla \Psi &=\nonumber  
                \Om\cdot \nabla^{\r} \Psi \\&+ \left[(1-\mu^2)\fbox{
  \begin{tabular}{c}
    \text{Is} $X^\n$ \\
    \text{flat?}
  \end{tabular}} + \mu\sqrt{1-\mu^2}\fbox{
  \begin{tabular}{c}
    \text{Is} $\gamma^\n$ \text{curved} \\
      \text{towards} $X^\n$?
  \end{tabular}} \right]\frac{\del \Psi}{\del \mu} \\&
\nonumber+ \left[ \sqrt{1-\mu^2}\left[\cosomega  \fbox{\begin{tabular}{c}
 \text{Is $\gamma^\t$}\\ \text{a geodesic} \\ \text{ of $\X^\n$?}
\end{tabular}} - \sinomega \fbox{\begin{tabular}{c}
 \text{Is $\gamma^\b$}\\ \text{a geodesic} \\ \text{ of $\X^\n$?}
\end{tabular}}\right] \nonumber  +  \mu \fbox{\begin{tabular}{c}
 \text{Does $\b$ not spin}\\ \text{around $\n$} \\ \text{ of $\X^\n$? }
\end{tabular}} \right]\frac{\del \Psi}{\del \omega}.  
        \end{align*}
    Since $\X^\n$ are spheres here the first answer to the first question is negative\footnote{{The last term is zero, if $\b$ spins around $\n$. The clumsy negation of the question is used to stay consistent with positive answers to these questions indicating the corresponding term being zero.}}. However, $\n$ is aligned with the radial vector everywhere, and hence $\gamma^\n$ is a straight line meaning the second term does indeed disappear. Next note that $\gamma^\t$ are in fact geodesics on the spheres, while $\b$ is constant along $\gamma^\n$. From this we conclude that all terms in $\nabla_{\Om}\mu$ except $(1-\mu^2)C(\omega)$ are zero as well as all terms in $\nabla_{\Om}\omega$ except $\sqrt{1-\mu^2}\sinomega (\t\cdot \kappab)$. \\
	
	Due to the symmetry of the sphere $C(\omega)= \frac{1}{\rho}$, however, is independent of $\omega$ and is just the curvature of the sphere. \\
	Considering the remaining term $\t\cdot \kappab$, note that $\gamma^\b$ at $\r$ are circles of radius $\sin(\theta)\rho$ with $\kappab$ pointing radially to the $z-$axis. Thus $|| \kappab||=\frac{1}{\sin(\theta)\rho}$ and therefore $\kappab\cdot\t=|| \kappab|| ||\t|| \cos(\theta) = \frac{\cot(\theta)}{\rho}$. 
	Piecing everything together 
	
	\begin{align}
		\nabla_{\Om}\Psi=& \nabla^\r_{\Om} \Psi + \frac{1-\mu^2}{\rho} \frac{\del \Psi}{\del \mu} 
		- \sqrt{1-\mu^2}\sin(\omega)\frac{\cot(\theta)}{\rho}\frac{\del f}{\del \omega}
	\end{align}

	\subsection{Ellipse}\label{sec:ex_ellipse}
	We will continue with an example that, on the one hand, is a natural next step from spherical coordinates, but on the other hand, will showcase a few practical issues that can arise as well as how quickly complexity can increase: Elliptical coordinates. For this example we will express equation \eqref{eq: OmPsi} as 
            \begin{align}
		\Om \cdot \nabla \Psi &=\nonumber  
                \Om\cdot \nabla^{\r} \Psi \\&+ \left[(1-\mu^2)\left[ \cos(\omega)^2 \n\cdot \kappat + \sin(\omega)^2 \n \cdot\kappab  \right]   + \mu\sqrt{1-\mu^2}\left( \hat{\Om}_\parallel\cdot\kappan\right) \right]\frac{\del \Psi}{\del \mu} \\&\nonumber+ \left[ \sqrt{1-\mu^2}\left[\cosomega(\b\cdot \kappat)- \sinomega (\t\cdot \kappab)\right] \nonumber  +  \mu (\t\cdot\nabla_{\n}\b) \right]\frac{\del \Psi}{\del \omega}.\label{eq:cyl_form}
        \end{align}
    We are about to see that even in this still rather simple case none of these terms will be zero for practical choices of $\n,\t$ and $\b$.
    
	Assume $\R^3$ is parametrized by $\rho, \phi, \theta$ again, such that
	$\r = \rho (a\cos(\phi)\sin(\theta),b \sin(\phi)\sin(\theta),c\cos(\theta))$, where $a,b,c>0$. A natural choice, and one that would get rid of $\del \Psi/\del \omega$ if ellipses are indeed level sets, would be to pick $\n$ as the normal vector to these ellipses. Thus \begin{align}
		\n \propto \left(\frac{\cosphi \sintheta}{a}, \frac{\sinphi \sintheta}{b},\frac{\costheta}{c}\right)^\top.
	\end{align}
	Easy choices for $\t$ and $\b$ would be 
	\begin{align}
		\t \propto \frac{\del \X}{\del \theta} = \left( a\cosphi\costheta, b\costheta\sinphi,c \sintheta \right)^\top
	\end{align}
	and \begin{align}
		\tilde{\b} \propto \frac{\del \X}{\del \phi} = \left(-a \sinphi\sintheta,b\cosphi \sintheta,0\right)^\top.
	\end{align}
	We will assume all three of these are normalized, even though we refrained from providing the normalization constants on the right-hand sides to simplify the presentation. This is essentially the same choice we made in the previous section on the sphere. In fact, if $a=b=c=1$ all vectors are identical. Note however that on the ellipsoid $\t$ and $\tilde{\b}$ are not always orthogonal unless $a=b$:
	\begin{align}
		&\t\cdot \tilde{\b} = \left[{(b^2-a^2)\cosphi\costheta\sinphi\sintheta}\right]\cdot\nonumber \\ &\left[\sin(\theta)^2\left(a^2\sinphi^2+b^2\cosphi^2\right) \left(a^2\cosphi^2\costheta^2 \right.\left. +b^2\costheta^2\sinphi^2+c^2\sintheta^2\right)  \right]^{-1/2}
	\end{align}
	
	The vector $\b = \n\times \t$ can be calculated as 
    \begin{align}
		\b & = \frac{\tilde{\b}-(\t\cdot\tilde{\b})\t}{1-\t\cdot \tilde{\b}} \label{eq:b_to_tilde} \propto \begin{pmatrix}-\frac{(b^2+c^2) \costheta\sinphi\sintheta}{bc} \\ \frac{(a^2+c^2)\cosphi\costheta\sintheta}{a c}\\-\frac{(a^2-b^2)\cosphi\cos(\theta)^2\sinphi}{ab}   \end{pmatrix}.
	\end{align}
    Note further that on an ellipsoid neither $\gamma^\t$ nor $\gamma^{\tilde{b}}$ are geodesics, forcing us to calculate the terms $\t\cdot\kappab$ and $\b\cdot\kappat$ or the corresponding shape operators.
	Although $\b$ can be used directly, many quantities, specifically $\gamma^\u$ and by extension $\kappa^{\u}$, are simpler to calculate for $\u= \tilde{\b}$ than for $\u = \b.$ Conveniently the map $\u \mapsto \kappa^{\u}$ is a linear mapping, as long as $\r$ is fixed, allowing us to work with $\tilde{\b}$, where useful, and then obtain $ \kappab$ via equation \eqref{eq:b_to_tilde}
	\begin{align}
		\kappab =  \frac{\bm{\kappa}^{\tilde{\b}}-(\t\cdot\tilde{\b}\kappat)}{1-\t\cdot\tilde{\b}}. \label{eq:kb_to_tilde}
	\end{align}
	The shape operator can be transformed in a similar fashion, if the surface curvature form is preferred.  The curvature $C(\omega)$ does not depend on the choice of basis of the tangent space. While it therefore does not depend on $\b$ or $\tilde{\b}$ it still implicitly depends on $\t$ through $\omega$. Because of this, it can be beneficial to use more convenient choices, such as $\tilde{\b}$, even if not all conditions are met, if appropriate adjustments are made.  We will demonstrate this in the current section. 
	
	In order to calculate $\kappan$, $\kappat$ and $\kappabt$ we need the curves $\gamma^\n$, $\gamma^\t$ and $\gamma^{\tilde{\b}}$.  Here is where chosing $\tilde{\b}$ over $\n\times \t$ becomes beneficial for now as parametrizations of $\gamma^\t$ and $\gamma^{\tilde{\b}}$ are given through $r(\rho,\phi,\theta)$ by fixing $\rho$ and one of the angles respectively.
	Thus $\gamma^\t(\tau,\r(\rho,\phi,\theta)) =  \r(\rho,\phi,\theta+\tau)$ and $\gamma^{\tilde{\b}}(\tau,\r(\rho,\phi,\theta)) =  \r(\rho,\phi+\tau,\theta)$ with $\tau\in (-\epsilon,\epsilon)$.  These are however not arc length parametrizations so either reparametrization is necessary or the formulas for ${\bm\kappa}^\u$ introduced in Remark \ref{rem:orth_frames} need adjustment. We correct for parametrization by setting \begin{align}
		{\bm \kappa}^\u = -\frac{(\gamma^\u)''}{|| (\gamma^\u)''||}\frac{||(\gamma^\u)'\times (\gamma^\u)'' ||}{|| (\gamma^\u)'||^3},
	\end{align}
	where the first term preserves the direction of the original definition, while the second term accounts for rescaling to arc-length.
	Using this we can calculate $\kappat$ and $\kappabt:$

	\begin{align}
	&	\kappat =\frac{1}{k_\t} \begin{pmatrix}
			{a \sin (\theta ) \cos (\phi ) \sqrt{c^2 \left(\left(a^2-b^2\right) \cos (2 \phi )+a^2+b^2\right)}}\\\\
			{b \sin (\theta ) \sin (\phi ) \sqrt{c^2 \left(\left(a^2-b^2\right) \cos (2 \phi )+a^2+b^2\right)}}\\\\
			{c \cos (\theta ) \sqrt{c^2 \left(\left(a^2-b^2\right) \cos (2 \phi )+a^2+b^2\right)}}
		\end{pmatrix}
	\end{align}
	with 
	\begin{align*}
	k_\t = \\ &\sqrt{2} \rho  \left(\cos ^2(\theta ) \left(a^2 \cos ^2(\phi )+b^2 \sin ^2(\phi )\right)+c^2 \sin ^2(\theta )\right)^{3/2}\sqrt{\sin ^2(\theta ) \left(a^2 \cos ^2(\phi )+b^2 \sin ^2(\phi )\right)+c^2 \cos ^2(\theta )}
	\end{align*}
	and

	\begin{align}
		&\kappabt = \frac{1}{k_{\tilde{\b}}} \begin{pmatrix}
			{a \sin (\theta ) \cos (\phi ) \sqrt{a^2 b^2 \sin ^4(\theta )}}
			\\\\{b \sin (\theta ) \sin (\phi ) \sqrt{a^2 b^2 \sin ^4(\theta )}}\\\\0
		\end{pmatrix},
	\end{align}
	where 
	\begin{align*}
		k_{\tilde{\b}} = \rho  \left(\sin ^2(\theta ) \left(a^2 \sin ^2(\phi )+b^2 \cos ^2(\phi )\right)\right)^{3/2}  \sqrt{\sin ^2(\theta ) \left(a^2 \cos ^2(\phi )+b^2 \sin ^2(\phi )\right)}
	\end{align*}
	Now $\kappab$ can be calculated through \eqref{eq:kb_to_tilde}, as can all inner products of $\kappat$ and $\kappab$ with $\n,\t$ and $\b$ required for the calculation $\nabla_{\Om}\Psi$.

	Finding a parametrization of $\gamma^\n$ to calculate $\kappan$ is less straightforward. Note in particular that fixing the two angles and varying $\rho$ does not give such a parametrization. A direct calculation of $\kappan = -\nabla_n\n$ can be simpler and is what we did to obtain
	
	\begin{align}
		\kappan= \frac{1}{k_\n} \begin{pmatrix}
			a \cosphi \sintheta \left(b^4 \left(c^2-a^2\right) \costheta^2\right.\left.+{c}^4 \left({b}^2-{a}^2\right) \sinphi^2 \sintheta^2\right),\\ \\
			b \sinphi \sintheta \left({a}^4 \left({c}^2-{b}^2\right) \costheta^2\right.\left.+{c}^4 \left({a}^2-{b}^2\right) \cosphi^2 \sintheta^2\right),\\\\
			{c} \sintheta^2 \costheta \left({a}^4 \left({b}^2-{c}^2\right) \sinphi^2\right.\left.-{b}^4 \left({c}^2-{a}^2\right) \cosphi^2\right)
		\end{pmatrix}
	\end{align}
	
	with \begin{align*} &k_\n =  {\rho} \left({c}^2 \sintheta^2 \left({a}^2 \sinphi+{b}^2 \cosphi^2\right)+{a}^2 {b}^2 \costheta^2\right)^2.\end{align*} 
	Calculating the inner product with $\Om_{\parallel}$ is thus simply a, somewhat tedious, exercise in arithmetic.
	Similarly, i.e. through direct calculation, we find the remaining term

	\begin{align}
		&\b\cdot \nabla_\n \t = -\t\cdot\nabla_\n \b =\nonumber\\
&\frac{{a b c \left(a^2-b^2\right) \sin (\phi) \cos (\phi) \cos (\theta)}}{\left[\rho \left(a^2 \sin(\phi)^2+b^2 \cos ^2(\phi)\right) \right] \left[a^2 b^2 \cos(\theta)^2 +c^2 \sin(\theta)^2 \left(a^2 \sin ^2(\phi)+b^2 \cos (\phi)^2\right)\right]}        
	\end{align}
	All that is left is to plug these terms back into the formulas for $\nabla_\Om\Psi$, which we will refrain from due to the length of the terms in this example. We saw that going from spherical to elliptical coordinates significantly increased the complexity of the involved terms as well as the procedure of calculating them. We also highlighted, that convenient choices might not always fulfill the conditions used in the derivations in the previous sections, but demonstrated how they might still work. 
	\begin{remark}{Homothetic Scaling}\
A shared feature of all the preceding examples is the appearance of the scaling factor $1/\rho$ in every term. This is not a coincidence, but rather a direct consequence of the fact that in each case, the rescaled embedding satisfies $X_{\rho \r}^\n = \rho X_\r^\n$. In other words, the foliation induced by $\n$ can be generated by scaling a single base manifold by a factor of $\rho$. Such families of manifolds are said to be generated by homothetic/self-similar rescaling. These commonly arise in coordinate systems and in the context of curvature flows. \cite{homothetical,Schrinkers}

More specifically, for $\r_1, \r_2 \in \mathbb{R}^3$ with $\r_2 = \rho \r_1$ for some $\rho>0$, we observe that the normal, tangent, and binormal vectors satisfy $\n(\r_1) = \n(\r_2), \t(\r_1) = \t(\r_2)$, and $\b(\r_1) = \b(\r_2)$. Consequently, the associated quantities $\X^\u$ and $\gamma^\u$ are simply rescaled by a factor of $\rho$ when transitioning from $\r_1$ to $\r_2$.

An illustrative 2D example of such homothetic rescaling is shown in Figure~\ref{fig:flower}.
		\begin{figure}
			\centering
			\includegraphics[width=.5\linewidth,trim={5cm 5cm 4cm 4cm},clip]{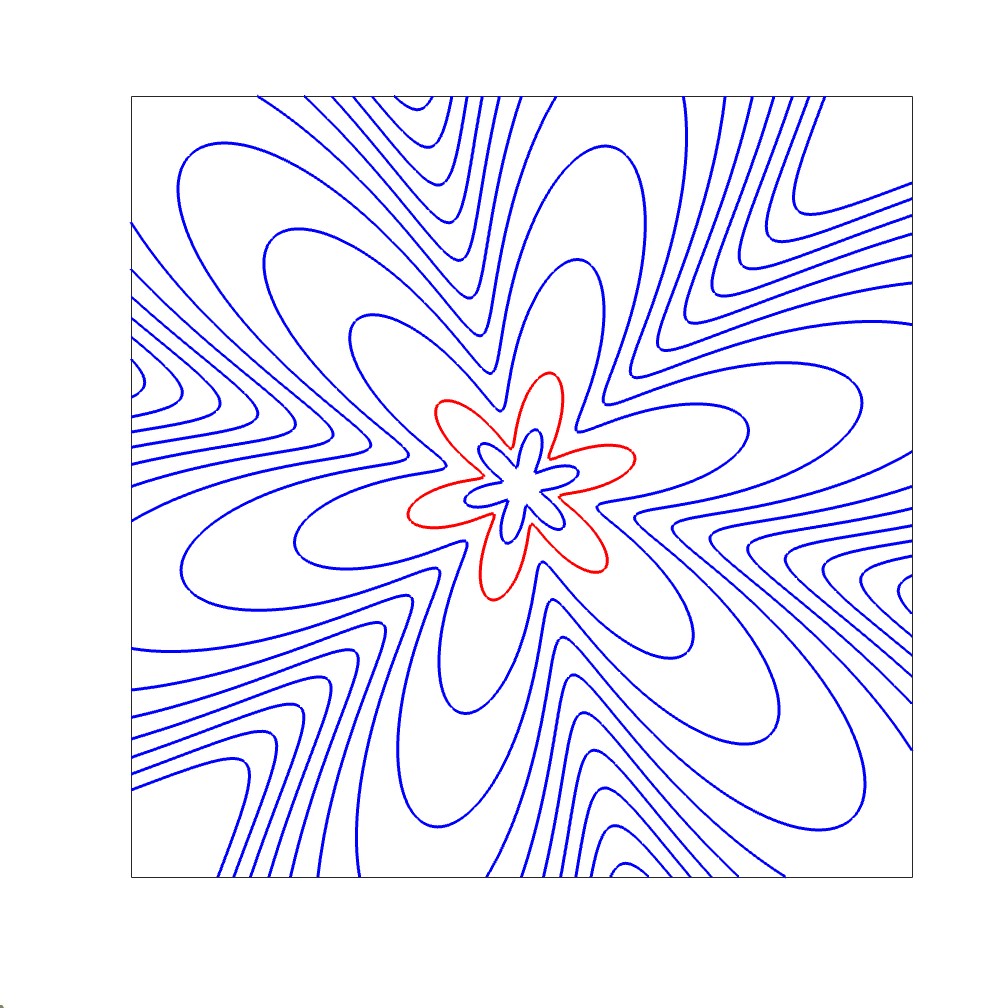}
			\caption{2D example of homothetically rescaled manifolds. Base manifold $\X$ in red, rescaled copies $\rho \X$ in blue for various $\rho.$}
			\label{fig:flower}
		\end{figure}
		
		To see the effect on curvatures let 
		$\X\in \R^3$ be an embedded manifold parametrized by $u^1,u^2 \in U^2\subset\R^2$, if $\X$ is a surface, and $u^1\in U^1\subset \R$, if $\X$ is a curve.
		\newcommand{\Y}{{\bm Y}}
		Next, let $\Y = \rho \X$. If $\X_i = \frac{\del \X}{\del u^i}$ with $i={1(,2)}$ are picked as tangent vectors to $\X$, then If $\Y_i = \frac{\del \Y}{\del u^i} = \rho \X_i$ are tangent vectors to $\Y$. 
		Since $\X_i$ and $\Y_i$ are colinear, the normal vectors $\N_\X=\N_Y$ are the same for fixed $u^1(,u^2)$. 
		Consequently  the metric tensors of $\Y$
		${g_\Y}_{\alpha,\beta} = \Y_\alpha \cdot \Y_\beta $
		$= \rho^2 \X_\alpha\cdot \X_\beta =$
		$ \rho^2{g_\X}_{\alpha,\beta}$ is a scalar multiple of the metric tensor $g_\X$ of $\X$. Similarly for the second fundamental form:
		$$ {h_\Y}_{\alpha,\beta} = \N_\Y\cdot \frac{\del^2 \Y}{\del \alpha \del \beta} = \rho N_\X\cdot \frac{\del^2 \X}{\del \alpha \del \beta} = \rho {h_\X}_{\alpha,\beta}.$$ This now finally leads to \begin{align}
			S_\Y = (g_\Y)^{-1}h_\Y = \frac{1}{\rho^2}(g_\X)^{-1} \rho h_\X = \frac{1}{\rho}S_\X.
		\end{align} 
		This relation between shape factors of scaled manifolds, intuitively implies that manifolds scaled by a factor $\rho$ are flatter by a factor $1/\rho$. This applies to manifolds in arbitrary dimensions, but specifically to the curves and surfaces of interest for us here thus explaining the factor $1/\rho$ in all terms from previous examples bar $\b\cdot \nabla_\n \t$.
		To explain the scaling of this factor consider $\gamma^\u(\r_1,\tau)$ and assume this is a parametrization by arc length. Then $\gamma^\u(\r_2,\tau) = \rho \gamma(r_1,\tau/ \rho)$ is also an arc length parametrization. Noting that this implies ${\bm \kappa}^\u(\r_1)= \frac{1}{\rho}{\bm \kappa^u(\r_2)}$ gives an additional explanation of the factor $1/\rho$ in all curvature related terms.
		However, we can also use this to see that 
		\begin{align}
			&\nonumber\nabla_\n\t(\r_2) \nonumber=\frac{d}{d\tau}\t(\gamma^n(\r_2,\tau))\bigg|_{\tau = 0} 
			=  \frac{d}{d\tau}\t(\rho\gamma^n(\r_1,\tau/\rho))\bigg|_{\tau = 0} \\ &= \frac{d}{d\tau}\t(\gamma^n(\r_1,\tau/\rho))\bigg|_{\tau = 0} =\frac{1}{\rho}\nabla_\n \t(\r_1)
		\end{align}  
		thereby explaining the $1/\rho$ factor in the remaining term. Here we used $\gamma^\u(\r_2,\tau) = \rho \gamma(r_1,\tau/ \rho)$ in the first line and $\t(\r) = \t(\rho \r)$ going from the first to the second line. 
	\end{remark}
	While the example of elliptical coordinates gave us a glimpse into how complex things can get one might find some relief in the fact that it suffices to handle the case $\rho=1$ and simply scale everything in the end, when dealing with coordinates based on homothetic foliations.

	\subsection{Translating Graph coordinates}\label{sec:ex_translators}
\begin{figure}
		\centering
		\includegraphics[width=.5\linewidth]{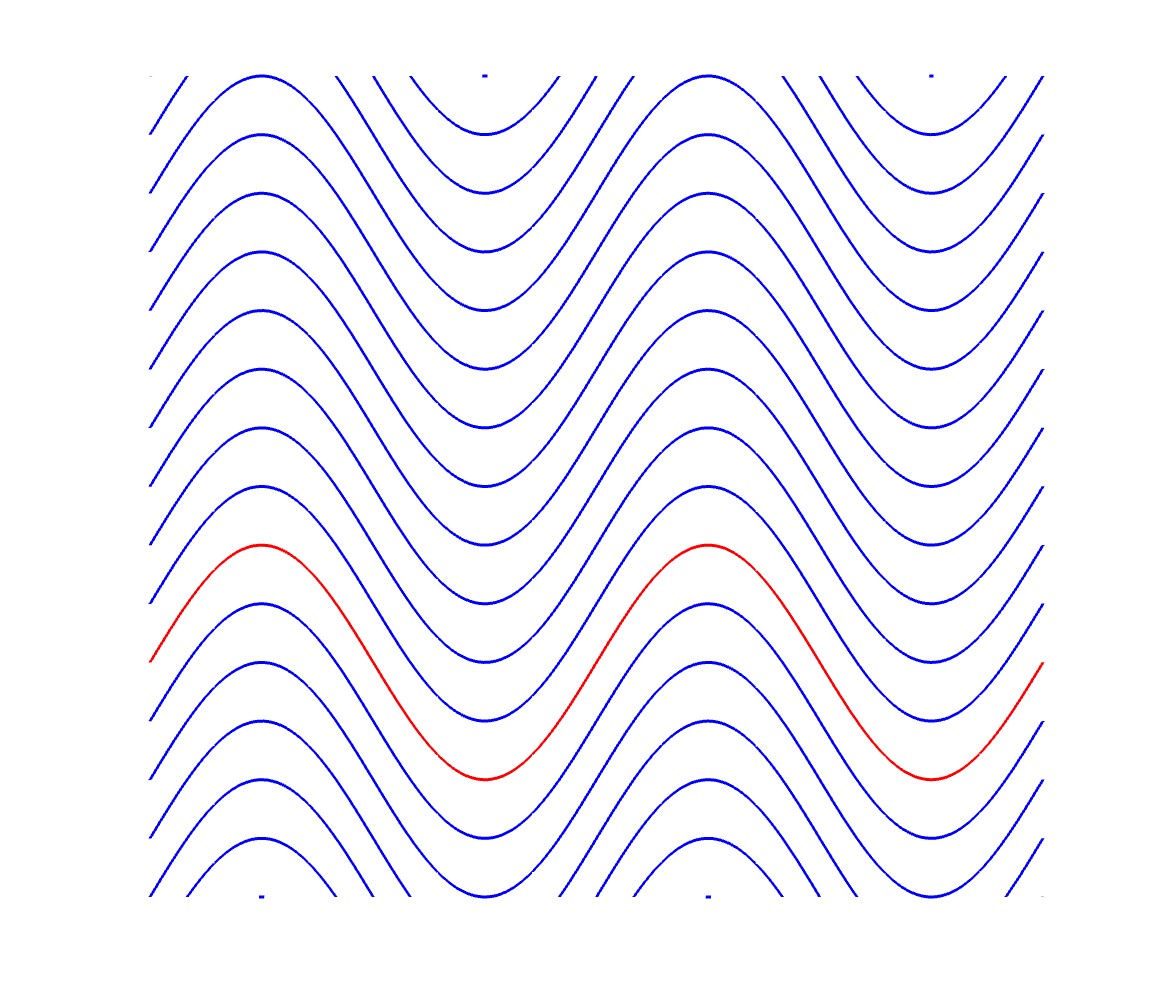}
		\caption{2D example of translating graphs. Base graph $G_0$ in red and several translated copies $G_z=G_0+z$ in blue.}
		\label{fig:translator}
	\end{figure}

	We will close the example section with a class of coordinates that is not based on homothetic foliation, but translated graphs instead. For this let $f\in C(\R^2,\R)$ be twice differentiable function with graph $G_0(x,y)=(x,y,f(x,y))$ and assume that $\R^3$ is parametrized along translated copies of $G$:  \begin{align}\r(x,y,z):=G_z(x,y):=(x,y,f(x,y)+z).\end{align} An example is depicted in figure \ref{fig:translator}. 
    We will consider the special case of a translating paraboloid first and then provide formulas for the case of general $G_z$. For these two examples we will use the form of $\Om\cdot\nabla\Psi$ given at the start of the example section in \eqref{eq:curv_form}. As was the case for ellipses unfortunately none of the terms will end up being zero. However we will see that all terms in the streaming operator can be calculated with knowledge of $G_0$ alone, i.e. independently of $z$. 
    \subsubsection*{Translating Paraboloid}
	In the case of $G_0$ being a paraboloid $f$ is of the form $f(x,y) = a x^2 + by^2$ ($a,b \in \R$). Tangent vectors parallel to the coordinate axis to $G_z$ can by found by normalizing $\del_x G_z$ and $\del_y G_z$ respectively.

This intuitive choice leads to
\begin{align}
    \t = \begin{pmatrix}
        \frac{1}{\sqrt{4 a^2 x^2+1}} \\0 \\\frac{2 a x}{\sqrt{4 a^2 x^2+1}}
    \end{pmatrix} &  & \tilde{\b}=\begin{pmatrix} 0 \\\frac{1}{\sqrt{4 b^2 y^2+1}} \\ \frac{2 b y}{\sqrt{4 b^2 y^2+1}}\end{pmatrix}
\end{align}
along with the normal vector 
\begin{align}
		\n(x,y,z)  
		 = \begin{pmatrix}
			-\frac{2 a x}{\sqrt{4 a^2 x^2+1} \sqrt{4 b^2 y^2+1}}\\-\frac{2 b y}{\sqrt{4 a^2 x^2+1} \sqrt{4 b^2 y^2+1}}\\\frac{1}{\sqrt{4 a^2 x^2+1} \sqrt{4 b^2 y^2+1}}
		\end{pmatrix}	
	\end{align}
   However, $\t$ and $\tilde{\b}$ are generally not orthogonal. This can be corrected via Gram-Schmidt, keeping $\t$ and replacing $\tilde{\b}$ with:
\begin{align}
    \b = \begin{pmatrix}
        -\frac{4 a b x y}{\sqrt{4 a^2 x^2+1} \left(\sqrt{4 a^2 x^2+1} \sqrt{4 b^2 y^2+1}-4 a b x y\right)}\\\frac{\sqrt{4 a^2 x^2+1}}{\sqrt{4 a^2 x^2+1} \sqrt{4 b^2 y^2+1}-4 a b x y}\\\frac{2 b y}{\sqrt{4 a^2 x^2+1} \left(\sqrt{4 a^2 x^2+1} \sqrt{4 b^2 y^2+1}-4 a b x y\right)}
    \end{pmatrix}
\end{align}
Using $\tilde{\b}$ over $\b$ has the advantage that $\kappabt = \del_y \tilde{\b}$, while $\kappab$ is a slightly more complicated directional derivative of $\b$. Similarly calculating the second fundamental form needed for the shape operator is slightly easier obtained with respect to $\{\t,\b\}.$
    The price to pay for this simplification is that terms calculated using $\tilde{\b}$ will require correction via  equations \eqref{eq:b_to_tilde} and \eqref{eq:kb_to_tilde}.

Either route leads to the shape operator with respect to $\t$ and $\b$ being

\begin{align}
   \left( S_{G_z}\right)_{\t\t} = & -\frac{2 a}{\left(4 a^2 x^2+1\right) \sqrt{4 b^2 y^2+1}} \\
   \left( S_{G_z}\right)_{\t\b} =& \frac{8 a^2 b x y \left(\sqrt{4 a^2 x^2+1} \sqrt{4 b^2 y^2+1}+4 a b x y\right)}{\left(4 a^2 x^2+1\right) \sqrt{4 b^2 y^2+1} \left(4 a^2 x^2+4 b^2 y^2+1\right)} \\
   \left( S_{G_z}\right)_{\b\t} = &\frac{8 a^2 b x y \left(\sqrt{4 a^2 x^2+1} \sqrt{4 b^2 y^2+1}-4 a b x y\right)}{\left(4 a^2 x^2+1\right) \sqrt{4 b^2 y^2+1} \left(4 a^2 x^2+4 b^2 y^2+1\right)} \\
    \left( S_{G_z}\right)_{\b\b} = &-\frac{2 b \left(\frac{16 a^3 b x^2 y^2}{\sqrt{4 b^2 y^2+1}}+4 a^2 x^2 \sqrt{4 a^2 x^2+1}+\sqrt{4 a^2 x^2+1}\right)}{\left(4 a^2 x^2+1\right) \left(4 a^2 x^2+4 b^2 y^2+1\right)}
   \end{align}
The curvatures $\kappan$,$\kappat$ and $\kappab$ are equally lengthy so that for the sake of readability we will only provide their inner products with $\Om_\parallel$,$\b$ and $\t$ as needed.
\begin{align}
    \Om_\parallel\cdot\kappan = &\cos(\omega)\frac{4 a^2 x}{\left(4 a^2 x^2+1\right)^{3/2} \left(4 b^2 y^2+1\right)} \nonumber\\
   & +\sin(\omega) \frac{4 b y \left(-4 a^3 x^2+4 a^2 b x^2+b\right)}{\left(4 a^2 x^2+1\right)^{3/2} \left(4 b^2 y^2+1\right) \left(\sqrt{4 a^2 x^2+1} \sqrt{4 b^2 y^2+1}-4 a b x y\right)}
\end{align}
Conveniently 
\begin{align}
		&\kappat\cdot \b = 0
	\end{align}
	and 
	\begin{align}
		\kappab\cdot \t =  \frac{4 a b y}{\left(4 a^2 x^2+1\right)^{3/2} \left(\sqrt{4 a^2 x^2+1} \sqrt{4 b^2 y^2+1}-4 a b x y\right)}
	\end{align}
	Lastly, we can find rather immediately 
	\begin{align}
		&	\b\cdot\nabla_\t \n = - \t \cdot \nabla_\b\n \nonumber=\\&
		\frac{8 a^2 b x y}{\left(4 a^2 x^2+1\right)^{3/2} \sqrt{4 b^2 y^2+1} \left(4 a b x y-\sqrt{4 a^2 x^2+1} \sqrt{4 b^2 y^2+1}\right)}
	\end{align}

    \subsubsection*{General translating Graph}
    As final example we shall generalize the previous example by replacing the paraboloid by a general function $f\in C^2(\R^2,\R)$. As shorthand let subscripts on $f$ denote derivatives with respect to the subscript, i.e. for example $f_x = \frac{\del f}{\del x}$. We will again chose $\n$ to be the normal vector of  the graph of $f$ given by
	 \begin{align}
		\n(x,y,z) = \begin{pmatrix}
			-\frac{f_x}{\sqrt{f_y^2+1} \sqrt{f_x^2+1}}\\ -\frac{f_y}{\sqrt{f_y^2+1} \sqrt{f_x^2+1}} \\ \frac{1}{\sqrt{f_y^2+1} \sqrt{f_x^2+1}}
		\end{pmatrix}.
	\end{align}
    Like for the paraboloid
	\begin{align}
		\t= \begin{pmatrix}
			\frac{1}{\sqrt{f_x^2+1}}\\0\\ \frac{f_x}{\sqrt{f_x^2+1}}
		\end{pmatrix} &  & \tilde{\b}=\begin{pmatrix}
			0\\\frac{1}{\sqrt{f_y^2+1}}\\\frac{f_y}{\sqrt{f_y^2+1}}
		\end{pmatrix}
	\end{align}
	are the canonical tangent vectors to $G_z$ parallel to the $x$ and $y$- axis and, also again, they generally will fail to be orthogonal.  
	An appropriate, orthonormal choice is \begin{align}
		\b(x,y) = \frac{1}{\sqrt{f_y^2+1} \sqrt{f_x^2+1}-f_y f_x}
		\begin{pmatrix}\frac{f_y f_x}{\sqrt{f_x^2+1} }\\\\ 
			-\sqrt{f_x^2+1}\\ \\
			\frac{f_y}{\sqrt{f_x^2+1} } \end{pmatrix}
	\end{align}
	
	In this basis we can find the components of the shape operator by calculating the metric tensor and second fundamental form and get 
	\begin{align}
		(S_{G_z})_{\t\t} &= \frac{-1}{k_S}f_{xx}\\ \nonumber\\
		(S_{G_z})_{\t\b}  &= \frac{f_y f_x f_{xx}-\sqrt{f_y^2+1} \sqrt{f_x^2+1} f_{xy}}{k_S \left(\sqrt{f_y^2+1} \sqrt{f_x^2+1}-f_y f_x\right)}\\\nonumber\\
        (S_{G_z})_{\b\t} & = 
        \frac{\left(\sqrt{f_y^2+1} \sqrt{f_x^2+1}-f_y f_x\right) \left(f_y f_x f_{xx}-\left(f_x^2+1\right) f_{xy}\right)}{k_S \left(f_y^2+f_x^2+1\right)}\\   \nonumber\\  
        (S_{G_z})_{\b\b} &= \frac{f_y f_x \left(\left(f_x^2+\sqrt{f_y^2+1} \sqrt{f_x^2+1}+1\right) f_{xy}\right. -f_y f_x f_{xx}-\sqrt{f_y^2+1} f_{yy} \left(f_x^2+1\right)^{3/2}}{
		k_S \left(f_y^2+f_x^2+1\right)}, 
	\end{align}
	where $$k_S = {\sqrt{f_y^2+1} \left(f_x^2+1\right)} .$$ 
	Again we will only provide the necessary inner products of $\kappan,\kappat$ and $\kappab$ with $\Om_{\parallel},\b$ and $\t$. 
	
	\begin{align}
		&\Om_{\parallel}\cdot\kappan = \cosomega \frac{f_y f_{xy}+f_x f_{xx}}{\left(f_y^2+1\right) \left(f_x^2+1\right)^{3/2}} +
		\\&\sinomega\nonumber \frac{\left(f_{yy} \left(f_x^2+1\right)-f_x^2 f_{xx}\right) f_y-f_y^2 f_x f_{xy}+f_x \left(f_x^2+1\right) f_{xy}}{\left(f_y^2+1\right) \left(f_x^2+1\right)^{3/2} \left(\sqrt{f_y^2+1} \sqrt{f_x^2+1}-f_y f_x\right)}
	\end{align}
	
	\begin{align}
		&\kappat\cdot \b  = -\frac{f_y \left(\sqrt{f_x^2+1} f_y^2+\sqrt{f_x^2+1}+\sqrt{f_y^2+1} \left(-f_y\right) f_x\right) f_{xy}}{\left(f_y^2+1\right) \left(f_x^2+1\right) \left(f_y f_x-\sqrt{f_y^2+1} \sqrt{f_x^2+1}\right)^2}
	\end{align}
	and 
	\begin{align}
		\kappab\cdot \t =  \frac{f_y f_{xx}}{\left(f_x^2+1\right)^{3/2} \left(\sqrt{f_y^2+1} \sqrt{f_x^2+1}-f_y f_x\right)}
	\end{align}
To wrap up we calculate the winding term as
    
    \begin{align}
		\b\cdot\nabla_\t \n = - \t \cdot \nabla_\b\n =
		-\frac{f_y \left(f_y f_{xy}+f_x f_{xx}\right)}{\sqrt{f_y^2+1} \left(f_x^2+1\right)^{3/2} \left(\sqrt{f_y^2+1} \sqrt{f_x^2+1}-f_y f_x\right)}
	\end{align}
    It should be noted that, while the general translating graph example expectedly looks more complicated than the paraboloid, the differences, formally, are not too big. Since the curvature terms consist of first and second derivatives only the paraboloid already exposes most of the complexity that is to be expected. Further note that none of the terms depend on $z$. Thus, despite none of the terms in these two example being particularly simple they all can be calculated using geometric properties of $G_0$ alone.
	
	\section{Conclusion}\label{eq:conclustion}
	We provided a derivation of the streaming operator $\nabla_{\Om}\Psi$ in equation \eqref{eq:kinetic_general} with a high level of generality and freedom concerning the choice of parametrization of $\Om$. We also provided various interpretations of the resulting terms relating to the curvature of manifolds, which we believe to be more intuitive than pure calculus or complicated pictures.  While we provided a handful of classical examples where the correctly parametrized form of $\nabla_{\Om}\Psi$ can be almost intuited, we also demonstrated the complexity can increase tremendously once the parametrization of $\Om$ becomes even slightly more complicated. Regardless, we believe the connection to curves and surfaces sheds light on the issue at hand from a new angle and might provide intuition and guidance to practitioners when deciding if the trade-off between increased complexity versus a potential reduction in dimension when using different coordinate systems is worth it or not. Beyond that we see use as a intuitive verification tool. While actual curvatures and thus the coefficients of terms in the streaming operator can be hard to calculate, whether or not they vanish can oftentimes be decided  without any calculation.
    
    \bibliographystyle{plainnat}
    \bibliography{references.bib}

\begin{thebibliography}{18}
\providecommand{\natexlab}[1]{#1}
\providecommand{\url}[1]{\texttt{#1}}
\expandafter\ifx\csname urlstyle\endcsname\relax
  \providecommand{\doi}[1]{doi: #1}\else
  \providecommand{\doi}{doi: \begingroup \urlstyle{rm}\Url}\fi

\bibitem[Arridge(1999)]{arridge1999}
Simon~R. Arridge.
\newblock Optical tomography in medical imaging.
\newblock \emph{Inverse Problems}, 15\penalty0 (2):\penalty0 R41--R93, 1999.

\bibitem[Bell and Glasstone(1970)]{Bell_Glasstone}
G.~I. Bell and S.~Glasstone.
\newblock \emph{Nuclear Reactor Theory}.
\newblock 1970, 10 1970.
\newblock URL \url{https://www.osti.gov/biblio/4074688}.

\bibitem[Bellomo et~al.(2013)Bellomo, Bellouquid, Tao, and Winkler]{bellomo2013}
Nicola Bellomo, Abdelghani Bellouquid, Yulong Tao, and Michael Winkler.
\newblock On the modeling of social crowds: A survey of models, speculations, and perspectives.
\newblock \emph{Mathematical Models and Methods in Applied Sciences}, 23\penalty0 (04):\penalty0 1230004, 2013.

\bibitem[Berger(1955)]{Berger}
Marcel Berger.
\newblock Sur les groupes d'holonomie homog\`ene des vari\'et\'es \`a{} connexion affine et des vari\'et\'es riemanniennes.
\newblock \emph{Bull. Soc. Math. France}, 83:\penalty0 279--330, 1955.
\newblock ISSN 0037-9484.
\newblock URL \url{http://www.numdam.org/item?id=BSMF_1955__83__279_0}.

\bibitem[Cercignani(1988)]{cercignani1988}
Carlo Cercignani.
\newblock \emph{The Boltzmann Equation and Its Applications}, volume~67 of \emph{Applied Mathematical Sciences}.
\newblock Springer, 1988.

\bibitem[Chandrasekhar(1960)]{chandrasekhar1960}
Subrahmanyan Chandrasekhar.
\newblock \emph{Radiative Transfer}.
\newblock Dover Publications, 1960.

\bibitem[Drugan et~al.(2018)Drugan, Lee, and Nguyen]{Schrinkers}
Gregory Drugan, Hojoo Lee, and Xuan~Hien Nguyen.
\newblock A survey of closed self-shrinkers with symmetry.
\newblock \emph{Results in Mathematics}, 73\penalty0 (1):\penalty0 32, 2018.
\newblock ISSN 1420-9012.
\newblock \doi{10.1007/s00025-018-0763-3}.
\newblock URL \url{https://doi.org/10.1007/s00025-018-0763-3}.

\bibitem[Freimanis(2011)]{FREIMANIS}
J.~Freimanis.
\newblock On vector radiative transfer equation in curvilinear coordinate systems.
\newblock \emph{Journal of Quantitative Spectroscopy and Radiative Transfer}, 112\penalty0 (13):\penalty0 2134--2148, 2011.
\newblock ISSN 0022-4073.
\newblock \doi{https://doi.org/10.1016/j.jqsrt.2011.04.007}.
\newblock URL \url{https://www.sciencedirect.com/science/article/pii/S0022407311001646}.
\newblock Polarimetric Detection, Characterization, and Remote Sensing.

\bibitem[Frobenius(1877)]{Frobenius}
G.~Frobenius.
\newblock Ueber das pfaffsche problem.
\newblock \emph{Journal für die reine und angewandte Mathematik}, 1877\penalty0 (82):\penalty0 230--315, 1877.
\newblock \doi{doi:10.1515/crll.1877.82.230}.
\newblock URL \url{https://doi.org/10.1515/crll.1877.82.230}.

\bibitem[Kohn(2008)]{homothetical}
Tobias Kohn.
\newblock Self-similar shrinkers in $\mathbb{R}^3$, January 2008.

\bibitem[Kreyszig(2013)]{kreyszig2013}
E.~Kreyszig.
\newblock \emph{Differential Geometry}.
\newblock Dover Books on Mathematics. Dover Publications, 2013.
\newblock ISBN 9780486318622.
\newblock URL \url{https://books.google.com/books?id=Vw3CAgAAQBAJ}.

\bibitem[Larsen(2007)]{helical}
Edward~W. Larsen.
\newblock The description of particle transport problems with helical symmetry.
\newblock \emph{Nuclear Science and Engineering}, 156\penalty0 (1):\penalty0 68--73, 2007.
\newblock \doi{10.13182/NSE07-A2685}.

\bibitem[Lee(2000)]{Lee00}
John~M. Lee.
\newblock \emph{Introduction to Smooth Manifolds}.
\newblock 2000.

\bibitem[Lee(2018)]{lee}
John~M. Lee.
\newblock \emph{Introduction to Riemannian Manifolds}.
\newblock Springer, Cham, Switzerland, second edition, 2018.
\newblock \doi{10.1007/978-3-319-91755-9}.
\newblock URL \url{https://doi.org/10.1007/978-3-319-91755-9}.

\bibitem[Majumdar(1993)]{majumdar1993}
Arun Majumdar.
\newblock Microscale heat conduction in dielectric thin films.
\newblock \emph{Journal of Heat Transfer}, 115\penalty0 (1):\penalty0 7--16, 1993.

\bibitem[Nicholson(1983)]{nicholson1983}
Dwight~R. Nicholson.
\newblock \emph{Introduction to Plasma Theory}.
\newblock John Wiley \& Sons, 1983.

\bibitem[Schlickeiser(2002)]{schlickeiser2002}
Reinhard Schlickeiser.
\newblock \emph{Cosmic Ray Astrophysics}.
\newblock Springer, 2002.

\bibitem[Weingarten(1861)]{Weingarten}
J.~Weingarten.
\newblock Ueber eine klasse auf einander abwickelbarer flächen.
\newblock \emph{Journal für die reine und angewandte Mathematik}, 59:\penalty0 382--393, 1861.
\newblock URL \url{http://eudml.org/doc/147844}.

\end{thebibliography}
	
	\appendix
	\section{Appendix}
	\subsection{When does $\del_\omega \Psi$ disappear?} \label{app:1}
	We mentioned in previous sections that if $\X^\n$ are picked to  be the level sets of $\Psi$, i.e. $\nabla^\r \Psi(\r,\Om,t)\propto \n(\r)$ 
	for all $\r,\Om$, then the terms $\frac{\del \Psi}{\del \omega}$ are also zero under rather mild conditions on $X^\n.$ These conditions are the following: Let $\v \in T_\r X^\n$ 
	be a tangent vector to $X^\n$ at $\r$ and let $Q\in SO(2)$. We require there exists a loop $\gamma:[0,1]\to X^\n$ on $X^\n$ starting and ending at $\r$, such that parallel transport \cite{lee} along $\gamma$ maps $v$ to $Qv$.  If this holds we say the Holonomy group  of $X^\n$ at $\r$ equals $SO(2)$, or $Hol_\r=SO(2)$.
	If this condition holds then for every $\omega_1$ and $\omega_2$, there is a loop $\gamma\subset X^\n$ such that $\Om(\mu,\omega_1)$ is turned into $\Om(\mu,\omega_2)$ while moving along $\gamma$. But since $X^\n$ was assumed to be a level set of $\Psi$ this implies $\Psi(\r,\Om(\mu,\omega_1),t)=\Psi(\r,\Om(\mu,\omega_2,t)$ and therefor $\frac{\del \Psi}{\del \omega}=0.$ While sounding rather complicated the condition $Hol_\r(\X^\n)=SO(2)$ holds for all orientable surfaces that are not a product space or have constant intrinsic curvature.\cite{Berger}
    While spheres are a notable exception whose holonomy group happens to be $SO(2)$ despite having constant curvature, cylinders unfortunately are not.
	
	\subsection{Some comments on Conservation forms}\label{app:2}
	In numerical applications people are often more interested in spatial and directional averages of $\Psi$ and as such of $\nabla_{\Om}\Psi.$ It is therefore convenient to rewrite $\nabla_{\Om}\Psi$ in conservation form \cite{Bell_Glasstone} 
	\begin{align}
		v(\r)\nabla_{\Om}\Psi = \begin{pmatrix}
			\nabla^\r \\ 
			\del_\mu \\
			\del_\omega
		\end{pmatrix}\cdot {\bm J}(\r,\mu,\omega), \label{eq:conservation}
	\end{align} 
	where $v(\r)$ is the volume element with respect to the coordinates of $\r$ and ${\bm J}\in \R^5$ is some vector field whose components can be interpreted as flux of $\Psi$. If such a form exists the integral 
	\begin{align}
		\int_V\int_{\S^2} \nabla_{\Om}\Psi = 	\int_V\int_{\S^2} \nabla_{\Om}\Psi v(\r) du^1du^2du^3 d\mu d\omega. 
	\end{align}
	is simplified tremendously via the divergence theorem. \\
	We show through a simplified case that a reformulation such as \eqref{eq:conservation} can not be expected in arbitrary coordinates.
	Suppose that we have picked coordinates $\r(u^1,u^2,u^3)$ such that $\frac{d \r}{du^1} = \n =\nabla \Psi$, that is $X^\n$ are level sets of $\Psi$. Further, assume that the structure of $\X^\n$ is sufficient for this to imply that $\del \Psi/\del \omega =0$. (Compare section \ref{app:1} for details. Also take note, that this example specifically excludes cylindrical coordinates.)
	
	We therefore have
	 \begin{align}
	&	v(\r)\nabla_\Om \Psi =v(\r)\mu \nabla_\n \Psi + v(\r)\left((1-\mu^2)C(\r,\omega)+\mu\sqrt{1-\mu^2}(\Om_\parallel\cdot\kappan)\right)\frac{\del \Psi}{\del \mu} \label{eq:cons1}
	\end{align}
	as shown in section 2. To arrive at a conservation form of this we attempt the following Ansatz:
	\begin{align}
	&	v(\r)\nabla_\Om \Psi  = \frac{v(\r)\mu}{f} \nabla_\n(f \Psi)  v(\r)\frac{\left((1-\mu^2)C(\r,\omega)+\mu\sqrt{1-\mu^2}(\Om_\parallel\cdot\kappan)\right)}{g}\frac{\del (g\Psi)}{\del \mu}\label{eq:cons2}
	\end{align}
	with some functions $f,g>0$. For this to be in conservation form three  conditions are required: \begin{align}
		&\nabla_\n\frac{v(\r)\mu}{f} = 0 \Rightarrow \nabla_\n \frac{v}{f}=0 \label{eq:cond1}\\ 
		&	\frac{\del }{\del \mu} \frac{\left((1-\mu^2)C(\r,\omega)+\mu\sqrt{1-\mu^2}(\Om_\parallel\cdot\kappan)\right)}{g} = 0  \label{eq:cond2}\\
		&	 \frac{\mu}{f} \nabla_\n f =  - \frac{\left((1-\mu^2)C(\r,\omega)+\mu\sqrt{1-\mu^2}(\Om_\parallel\cdot\kappan)\right)}{g}\frac{\del g}{\del \mu} \label{eq:cond3}
	\end{align}
	Where the last condition merely ensures equality of the right-hand sides of \eqref{eq:cons1} and \eqref{eq:cons2}. 
	
	Through some manipulations and after abbreviating $$K(\r,\mu,\omega)=\left((1-\mu^2)C(\r,\omega)+\mu\sqrt{1-\mu^2}(\Om_\parallel\cdot\kappan)\right)$$ conditions \eqref{eq:cond1} and \eqref{eq:cond2} can be rewritten as 
	\begin{align}
		\frac{\nabla_\n v}{v} = \frac{\nabla_\n f}{f}  \label{eq:cond1b}\\
		\frac{\del K}{\del \mu} = \frac{K}{g}\frac{\del g}{\del \mu} \label{eq:cond2b}
	\end{align}
	Plugging \eqref{eq:cond2b} into \eqref{eq:cond3} then leads to \begin{align}
	&	\frac{\nabla_\n v}{v}=\frac{\nabla_\n f}{f} = \frac{1}{\mu}\frac{\del K}{\del \mu}  \nonumber = C(\r,\omega)-\frac{1-2\mu^2}{\mu\sqrt{1-\mu^2}}(\Om_{\parallel}\cdot \kappan). \label{eq:cond_final}
	\end{align}
	In this equation the left-hand side is independent of $\mu$ and $\omega$, implying it can only hold if $\kappan = 0$ for all $\r$ and if $C(\r,\omega)=C(\r)$ is itself independent of $\omega$. 
	The latter condition means the sectional curvature of $\X^\n$ is the same in all directions, which from embedded surfaces leaves only planes and spheres. In both cases, conservation forms are known:
	\begin{align}
		f_{\text{plane}}=1 & &g_{\text{plane}}=1 \\
		f_{\text{sphere}} = \rho & &g_{\text{sphere}} = 1.
	\end{align}
	All other cases fail to produce a conservation form via the chosen assumptions and Ansatz. Note that, while cylindrical coordinates,  technically do not fulfill all assumptions used here, there still are conservation forms in the literature. These are however also based on a choice of $\n$ where $X^\n$ are planes. 
Lastly note that this is not proof that conservation forms do not exist in general, but merely that, if they do, finding them is more complicated than using a simple Ansatz like we did. 
\end{document}